\input amstex      
 \documentstyle{amsppt}
 \magnification 1200 

\hoffset= -.05 truecm
\hsize= 16.5 truecm
\vsize= 22 truecm
\voffset= -.05 truecm

 \def\ns{\hskip -.15 truecm} 
 \NoBlackBoxes

\def\summary{Theorem }
 
 \def\notation{Notation }
 \def\defcan{Definition 1.1 } 
\def\remcan{Remark 1.2 }
 \def\split{Proposition 1.3 } 
 \def\regPtwo{Proposition 1.4 }
 \def\hypsplit{Proposition 1.5 }   
 \def\splitwist{Proposition 1.6 } 
 \def\NoVer{Corollary 1.7 }  
  
 \def\defGal{Definition 2.1 }  
 \def\Galremtwo{Remark 2.2 } 
 \def\Galemma{Lemma 2.3 }
 \def\Galquad{Proposition 2.4 }
 \def\Galquadone{(2.4.1) }
 \def\defnrem{Definition-Remark 2.5 } 
 \def\bipro{Proposition 2.6 }
 \def\Galcyclic{Proposition 2.7 }
 \def\Notarem{Notation-Remark 2.8 } 
 \def\biproconv{Proposition 2.9 }
 \def\Galcycliconv{Proposition 2.10 } 
 \def\strPtwo{Theorem 2.11 }
 \def\Ptwosing{Corollary 2.12 }
 \def\Ptwoexamples{Proposition 2.13 }

 \def\scrbid{Theorem 3.1 }  
 \def\propnoprime{Theorem 3.2 }
\def\cornoprime{Corollary 3.3 }
 \def\ques{Question 3.4 }
 \def\nosimcyclreg{Theorem 3.5 }
 \def\nosimcyclone{(3.5.1) }
 \def\nosimcyclirreg{Theorem 3.6 }
 \def\nosimcycltwo{(3.6.1) }
 \def\nosimcyclsing{Corollary  3.7 }
 \def\strsmscreg{Theorem 3.8 }
 \def\strsmscrirr{Theorem 3.9 }
 \def\strsmscrirrnumber{(3.9.1) }
\def\classmooth{Theorems 2.11, 3.5, 3.6, 3.8 and 3.9 }
\def\classcroll{Theorems 3.5, 3.6, 3.8 and 3.9 }
 \def\splitscr{Corollary 3.10 }
 \def\exampcyc{Proposition 3.11 }
 \def\exampbid{Proposition 3.12 }

\def\Be{[Be] }
 \def\BS{[BS] }
 \def\Ca{[C] }
 \def\GPcy{[GP1] }
 \def\GPtrans{[GP2] }
\def\CR{[GP3] }
\def\GPsing{[GP4] }
 \def\GPring{[GP5] }
 \def\Ga{[Ga] }
 \def\Gr{[Gr] }
 \def\Hoone{[Ho1] }
 \def\Hotwo{[Ho2] }
 \def\Hothree{[Ho3] }
 \def\Hofour{[Ho4] }
 \def\HM{[HM] }
 \def\Kon{[K] }
 \def\MP{[MP] }
 \def\OP{[OP] }
 \def\Pa{[Pa] }
 \def\Pe{[Pe] }
 \def\Pu{[Pu] }
 \def\T{[T] }

 \topmatter
 \title Classification of quadruple Galois canonical covers I
 \endtitle
 
 \author Francisco Javier Gallego \\ and \\ Bangere  P. Purnaprajna
 \endauthor
 \abstract{In this article we classify quadruple Galois 
 canonical covers of smooth surfaces of minimal degree. The classification 
 shows that they are either non-simple cyclic covers or bi-double covers.
 If they are bi-double then they are all fiber products 
 of double covers. We construct examples to show 
 that all the possibilities in the classification do exist. There are 
 implications of this classification that include the existence of 
 families with unbounded geometric genus  and families with unbounded
 irregularity, in  
 sharp contrast with the case of double and triple canonical covers. 
 Together with the
results of Horikawa and Konno for double and triple 
 covers, a pattern emerges that motivates some
 general questions on the existence of higher degree canonical covers,
 some of which are answered in this article.}
 \endabstract
  \endthanks
 \address{Francisco Javier Gallego: Dpto. de \'Algebra,
  Facultad de Matem\'aticas,
  Universidad Complutense de Madrid, 28040 Madrid,
 Spain}\endaddress
 \email{gallego\@mat.ucm.es}\endemail
 \address{ B.P.Purnaprajna:
 405 Snow Hall,
   Dept. of Mathematics,
   University of Kansas,
   Lawrence, Kansas 66045-2142}\endaddress
 \email{purna\@math.ukans.edu}\endemail
 \thanks{The first author was partially supported by MCT project
number BFM2000-0621. He  is grateful for the hospitality of the
Department of Mathematics of the University of Kansas at
Lawrence. The second author is grateful to NSA for supporting this
research project. He is also grateful for the hospitality of the
Departamento de \'Algebra  of the Universidad Complutense de Madrid}
\endthanks

 \endtopmatter
 \document
 
 \vskip .2 cm
 
 \headline={\ifodd\pageno\rightheadline \else\leftheadline\fi}
 \def\rightheadline{\tenrm\hfil \eightpoint CLASSIFICATION OF QUADRUPLE
   GALOIS CANONICAL COVERS I
  \hfil\folio}
 \def\leftheadline{\tenrm\folio\hfil \eightpoint F.J. GALLEGO \&
   B.P. PURNAPRAJNA \hfil}
 
 \heading Dedicated to Ignacio Sols. \endheading
 
 \vskip .3 cm
 
 \heading Introduction \endheading

Classification problems are of central importance in Algebraic Geometry.
In the realm of algebraic surfaces, the geography of surfaces of general type,
by far the largest class of surfaces, is much less charted and understood. 
An important sub-class of surfaces of general type are those whose canonical 
map is a cover of a simpler surface, most notably, a cover of a surface of 
minimal degree. 
In the seventies and eighties Horikawa and Konno (see \Hoone \ns, 
\Hothree and \Kon \ns)
classified these covers when the degree of the cover 
is $2$ and $3$. In this article and in its sequel \GPsing  
we classify surfaces of general type whose canonical map is a 
 quadruple Galois cover of a surface of minimal degree. 

\smallskip

Covers of varieties of minimal degree have a ubiquitous presence in 
 various contexts. 
They appear 
 in the classification of surfaces of general type $X$ with small $K_X^2$
 done by 
 Horikawa (see \Hoone \ns , \Hotwo \ns, \Hothree and \Hofour \ns) and 
play an important role in mapping the
geography of surfaces of  
 general type.
They are also the chief source in constructing new examples of 
 surfaces of general type as the work of various geometers illustrates
 (see also \Kon \ns, \MP and \T \ns). 
These covers occur as well 
in the study of linear series on important 
 threefolds such as Calabi-Yau threefolds as the work in \BS \ns, \GPcy 
 and \OP shows.  
 They also become relevant in the study of the canonical ring of a variety 
 of general type as can be seen in results from \GPtrans and \Gr  
 (see \Pu for further motivation). 
 
 \smallskip

Compared to the canonical morphism of a curve, the
canonical morphism of a surface is much more subtle and allows a much
wider range of possibilities due to the existence of higher degree covers. 
The degree of the canonical morphism of a curve is bounded by $2$.   
In contrast, Beauville proved that 
the degree of the canonical morphism from a surface of
general type $X$  
onto a surface of minimal degree, or more generally, onto a surface 
 with geometric genus $p_g=0$, is bounded by $9$ if $\chi (X) \geq
 31$ (see \Be \ns; see also \Ga \ns, where the bound is reduced to $8$ 
if $p_g (X) \geq 132$).

\medskip

Surfaces of minimal degree are classically known to be linear $\bold P^2$, 
the Veronese surface in $\bold P^5$ and rational normal scrolls, which can 
be smooth (these include the smooth quadric hypersurface in $\bold P^3$) 
or singular (these are cones over a rational normal curve). As pointed out 
before the classification of the canonical covers of these surfaces is 
only complete when the degree of the cover is $2$ (see \Hoone \ns) and 
$3$ (see \Hothree and \Kon \ns). Horikawa also studied quadruple covers of 
linear $\bold P^2$. The next step in this classification is the study 
of quadruple covers of an arbitrary surface of minimal degree. 

\medskip

In this work we classify all quadruple Galois canonical 
covers of smooth surfaces of minimal degree $W$. In \GPsing we  
classify quadruple Galois covers of $W$ when $W$ is singular. 
There are many interesting consequences of the classification done here 
and in \GPsing \ns. Our classification yields, among other things, some 
striking contrasts with double and triple covers. 
Before we look at them we state the main result of this 
article: 

\proclaim{\summary }  Let $X$ be a canonical surface and let 
$W$ be a smooth surface of minimal degree. 
If the canonical bundle of $X$ is base-point-free and $X @> \varphi >> W$ is  
a quadruple Galois canonical cover,
then $W$ is either linear $\bold P^2$ or a smooth rational normal
scroll $S(m-e,m)$ with $0 \leq e \leq 2$ and $m \geq e+1$. Let
$G$ be the Galois group of $\varphi$.  

\smallskip

\item{A)} If $G= \bold Z_2 \times \bold Z_2$, then 
$X$ is the fiber product over $W$ of two 
double covers of $W$ branched along divisors $D_1$ and $D_2$ and $\varphi$ 
is the natural morphism from the fiber product to $W$. The divisors
$D_1$ and $D_2$ satisfy:

\smallskip

\item{1)} If $W$ is linear $\bold P^2$, then $D_1$ and $D_2$ are
  quartic curves. 

\item{2)} If $W=S(m-e,m)$, then either

\smallskip

\item{}2.1) $D_1 \sim 2C_0+(2m+2)f$ and $D_2 \sim 4C_0+(2e+2)f$;
or

\smallskip

\itemitem{} $e=0$ and 

\smallskip

\item{}2.2) $D_1 \sim (2m+2)f$ and $D_2 \sim 6C_0+2f$; or

\item{}2.3) $D_1 \sim (2m+4)f$ and $D_2 \sim 6C_0$; or

\item{}2.4) $D_1 \sim 2C_0+(2m+4)f$ and $D_2 \sim 4C_0$. 

\medskip

\item{B)} If $G=\bold Z_4$, then 
$\varphi$ is the composition of two double covers 
$X_1 @> p_1 >> Y$ branched along a divisor $D_2$
and $X @> p_2 >> X_1$, branched along the 
ramification of $p_1$ and $p_1^*D_1$ and with trace zero module
$p_1^*\Cal O_W(-\frac 1 2 D_1-\frac 1 4 D_2)$.  
The divisors $D_1$ and $D_2$ satisfy:

\smallskip

\item{1)} If $W$ is linear $\bold P^2$, $D_1$ is a conic  and $D_2$ is
  a quartic curve.

\item{2)} If $W=S(m-e, m)$, then either

\smallskip

\item{}2.1) $D_1 \sim (2m-e+1)f$ and $D_2 \sim 4C_0+(2e+2)f$; or 

\smallskip

\itemitem{} $e=0$, and

\smallskip

\item{}2.2) 
$D_1 \sim (2m+4)f$, $D_2 \sim 4C_0$; or 

\item{}2.3) $D_1 \sim 3C_0$, $D_2 \sim 2C_0+4f$ and $m=1$.

\medskip

\noindent Conversely, if 
$X @> \varphi >> W$  is either 

\smallskip

\item{I.} the fiber product over $W$ of two 
double covers $X_1 @> p_1 >> W$ and $X_2 @> p_2 >> W$, branched
respectively along divisors $D_2$ and $D_1$ as described 
in A.1, A.2.1, A.2.2, A.2.3 or A.2.4; or 

\item{II.}  the composition of two double covers 
$X_1 @> p_1 >> Y$ branched along a divisor $D_2$
and $X @> p_2 >> X_1$, branched along the 
ramification of $p_1$ and $p_1^*D_1$ and with trace zero module
$p_1^*\Cal O_W(-\frac 1 2 D_1-\frac 1 4 D_2)$, with $D_1$ and $D_2$ as 
described in B.1, B.2.1, B.2.2 or B.2.3, 

\smallskip

\noindent then $X @> \varphi >> W$ is a Galois canonical cover whose Galois
group is $\bold Z_2 \times \bold Z_2$ in case I and $\bold Z_4$ in
case II. 
\endproclaim

The above result is split inside the paper into  \NoVer
and \classmooth \ns.  
We also construct families of examples to show the existence of all the 
cases that appear in the classification.
A notable fact that is not mentioned in the statement of the theorem
is that there exist families of smooth surfaces of general type $X$ for 
every case A.1, A.2.1, A.2.2, A.2.3 and A.2.4, (that is, all cases where 
$G=\bold Z_2 \times \bold Z_2$). In comparison, we show that the quadruple 
cyclic canonical covers of smooth surfaces of minimal degree are always 
singular. We also show that quadruple cyclic canonical covers are 
non-simple cyclic.

\smallskip

One of the important implications of the main result of this article
is the existence of families of quadruple canonical covers with 
unbounded geometric genus and the existence of families with unbounded
irregularity. 
The unboundedness of the geometric genus is in sharp contrast
with the situation of triple covers and the unboundedness of the
irregularity is in sharp contrast with the situation of both double 
and triple covers.
The geometric genus of canonical double covers is unbounded but they are 
all regular, and even simply connected surfaces. The geometric
genus of canonical triple covers is bounded by $5$ (see \Kon \ns) and
are all regular. For quadruple Galois covers we show    
the existence of families of
surfaces $X$ for each possible value of $p_g(X)$.  
Regarding irregularity, the classification shows that most of the quadruple 
Galois covers are regular. However, we do show that there exist families of 
surfaces with irregularity $1$ (cases A.2.4, B.2.2 and B.2.3) and, most 
importantly, families of surfaces $X$ for each possible value of $q(X)$ 
(cases A.2.2 and A.2.3).

\smallskip

The classification of quadruple covers provides other significant contrast 
with double and triple covers and brings out clearly the marked difference
between even and odd degree covers. The only smooth targets of quadruple 
Galois canonical covers that occur are linear $\bold P^2$ and rational 
normal scrolls which correspond to only three Hirzebruch surfaces, namely 
$\bold F_0$, $\bold F_1$ and $\bold F_2$.
In the case of canonical double covers, linear 
$\bold P^2$ and smooth rational scrolls corresponding to every 
Hirzebruch surface appear as image of the canonical morphism. In the 
case of canonical triple covers, the list is reduced drastically and 
the only possible smooth target is linear $\bold P^2$ (see also 
\GPtrans \ns, Proposition 3.3).  
\smallskip

The classification of quadruple covers of singular targets (i.e., a cone 
over a rational normal curve) in \GPsing together with the results of
Horikawa and Konno exhibit a striking pattern.  
Indeed we show that these form a bounded family with respect to both the
geometric genus and irregularity as in the case of double and triples covers
of singular targets.
The results in this article and \GPsing predict a precise numerology 
that might hold for higher degree covers. 
The following facts make it clear what we mean:
there do not exist canonical covers of odd degree 
of smooth scrolls 
(see \text{\GPtrans \ns,} Proposition 3.3) and there do not exist 
Galois canonical covers $X @> \varphi >> W$ of prime degree $p > 3$ of 
surfaces $W$ of minimal degree if $X$ is regular or if $W$ is smooth 
(see \propnoprime \ns, \cornoprime and \GPtrans \ns, Corollary 3.2). 
This motivates us to pose a general question (see \ques \ns) 
on the non-existence of higher, prime degree Galois canonical covers of 
surfaces of minimal degree.

\smallskip

The classification obtained in this article and in \GPsing has 
further applications. In \GPring \ns, we determine the ring generators 
of the quadruple covers classified here  and in \GPsing \ns. 
The results of \GPring 
show that quadruple covers serve as examples and counter 
examples to some questions on graded rings and normal generation of 
linear systems on an algebraic surface.

 \medskip

 {\bf Acknowledgements:} We thank N. Mohan Kumar for generously sharing his 
time for some useful 
discussions and insightful comments. 
We are very grateful to Rita Pardini, who, 
after seeing an earlier version of this work, remarked that
we had not treated the case of non-simple cyclic covers. It was after
her  
comment that we proved the theorems on cyclic covers.
Finally, we thank the referee of the announcement \CR  of these results, 
who kindly made many suggestions, not only for the announcement, but 
also to improve the exposition of the present article and \GPsing \ns.

 \heading 1. Quadruple canonical covers of surfaces of minimal degree
 \endheading 
 
 \noindent {\bf Convention:} 
 We work over an algebraically closed field of
 charactheristic $0$. 

\medskip
 
 \noindent{\bf{\notation}}\ns : 
  {\it We will follow these conventions:} 
 
 \smallskip
 
\noindent  1) Throughout this article, unless otherwise stated,  
$W$ will be an embedded smooth projective
   algebraic surface of 
 minimal degree, i.e., whose degree is equal to its codimension in
 projective space plus $1$. 
 
 \smallskip
 
\noindent  2) Throughout
 this article,  unless otherwise stated, $X$ will be a projective
 algebraic normal 
 surface 
with at worst canonical singularities. 
 We will denote by  $\omega_X$ the canonical bundle of $X$.

 \smallskip
 
 \noindent We recall the following standard notation: 
 
\smallskip

\noindent  3) By $\bold F_e$ we denote the Hirzebruch surface whose
 minimal section have self-intersection $-e$. If $e >0$ let $C_0$
 denote the minimal section of $\bold F_e$ and let $f$ be one of the
 fibers of $\bold F_e$. If $e=0$, $C_0$ will be a fiber of one of the
 families of lines and $f$ will be a fiber of the other family of
 lines. 
 
 \smallskip
 
\noindent  4) If $a, b$ are integers such that $0 < a \leq b$, consider
 two disjoint linear subspaces $\bold P^a$ and $\bold P^b$ of $\bold
 P^{a+b+1}$. We denote
   by $S(a,b)$ the smooth rational normal scroll obtained by joining
   corresponding points of a rational normal curve in $\bold P^a$
   and a rational normal curve of $\bold P^b$. 
Recall that $S(a,b)$ is the
image of $\bold F_e$ by the embedding induced by the complete linear series
$|C_0+mf|$, with $a=m-e$, $b=m$ and $m\geq e+1$. 
If $a=b$, the linear series 
$|mC_0+f|$ also gives a minimal degree embedding of $\bold F_0$, 
equivalent to the previous one by
the automorphism of $\bold P^1 \times \bold P^1=\bold F_0$ swapping the
factors. In this case {\it our convention will always be to choose $C_0$ and
$f$ so that, when $W$ is a smooth rational normal scroll, and $W$ is
embedded by $|C_0+mf|$.}

 \medskip

 \medskip
 
 In this section we prove general results regarding quadruple canonical
 covers of surfaces of minimal degree: 
 
 \proclaim{\defcan} Let $X$ and $W$ be as in the previous notations. We
 will say that a surjective morphism $X @> \varphi >> W$ is a canonical
 cover of $W$ if $X$ is surface of general type whose canonical
 bundle $\omega_X$ is ample  
and
base-point-free and
 $\varphi$ is the canonical morphism of $X$. 
 \endproclaim

 \proclaim{\remcan} Although we have assumed $X$ to have canonical
 singularities, 
some results hold in greater generality. Precisely,
 if for the purpose of this remark we ignore Notation 2) above and $X$ is
 assumed to be a normal, locally Gorenstein surface 
 instead
 then \defcan still makes sense and 
 \classmooth  hold. We can further relax the hypotheses on $X$ in 
 the converse parts of \classmooth and they hold if $X$ is just assumed to be
 smooth in codimension $1$, since in that case 
 these covers are Gorenstein. 
 \endproclaim

 \proclaim{\split} Let $X @> \varphi >> W$ be a quadruple canonical 
cover of $W$. Let
 $H=\Cal O_W(1)$. 
Then $\varphi_*\Cal O_{X}$ is a vector bundle on $W$ and 
 $$\varphi_*\Cal O_X=\Cal O_W \oplus E \oplus (\omega_W
 \otimes H^*)$$
  with $E$ vector bundle over $W$ of rank $2$. 
 If in addition $\varphi_*\Cal O_{X}$ splits as sum of line
 bundles, then 
 $$\varphi_*\Cal O_{X}=\Cal O_W \oplus L_1^* \oplus L_2^* \oplus (\omega_W
 \otimes H^*) $$
 with $L_1^* \otimes L_2^*= \omega_W
 \otimes H^*$.  
 \endproclaim
 
 \noindent {\it Proof.} Recall that by 
 \defcan $\varphi$ is finite, $W$ is smooth and 
 $X$ is locally Cohen-Macaulay. Then  $\varphi$ is flat and 
 hence $\varphi_*\Cal O_{X}$ is a vector bundle over $\Cal O_W$ of rank
 $4$. Moreover (see \HM \ns), $\varphi_*\Cal O_{X}=\Cal O_W \oplus E'$, where $E'$
 is the trace zero module of $\varphi$. From relative duality we have
 $$\varphi_*\omega_{X}=(\varphi_*\Cal O_{X})^* \otimes 
              \omega_W \ .$$ On the other hand, by hypothesis,
              $\omega_{X}=\varphi^*H$, hence, by projection formula,  
 $$\varphi_*\omega_{X}=\varphi_*\Cal O_{X} \otimes H \ .$$
 Then $$\omega_W \oplus (\omega_W \otimes (E')^*)= 
 H \oplus (E' \otimes H) \ .$$
 Since $\omega_W=H$ is not possible, for $W$ is a rational surface,
 $E'=E \oplus (\omega_W \otimes H^*)$, with $E$ vector bundle of rank
 $2$. If $\varphi_*\Cal O_X$ splits, let $E=L_1^* \oplus L_2^*$.   
 Then 
  $$\displaylines {\omega_W \oplus (\omega_W \otimes L_1) \oplus
 (\omega_W \otimes L_2) \oplus 
 H= \cr
 H \oplus (H \otimes L_1^*) \oplus (H \otimes L_2^*) \oplus \omega_W   \ .}$$
Then taking the determinant of both sides of the equality gives
$(L_1^* \otimes L_2^*)^{\otimes 2}=(\omega_W \otimes H^*)^{\otimes
  2}$. Since $W$ is either $\bold P^2$ or a Hirzebruch surface, then 
$L_1^* \otimes L_2^* = \omega_W \otimes
 H^*$. \qed

 \bigskip
 
 Now we study in more detail the possible splittings of $\varphi_*\Cal O_X$
 depending on what surface $W$ is. We start with this observation about 
 linear $\bold P^2$:
 
 \proclaim{\regPtwo } Let $X @> \varphi >> W$ be a 
canonical
 cover. 
 If $W$ is linear $\bold P^2$, then $X$
 is regular if and only if $\varphi_*\Cal O_W$ splits as direct sum 
 of line bundles.
 \endproclaim
 
 {\it Proof.} By \split we know that
  $$\varphi_*\Cal O_X=\Cal O_W \oplus E \oplus \omega_W(-1) \ .$$
 Since $W=\bold P^2$, the intermediate cohomology of line bundles on
 $W$ vanishes so by projection formula 
 $H^1(\varphi^*\Cal O_W(k))=H^1(E(k))$. By Kodaira vanishing and
 duality $H^1(\varphi^*\Cal O_W(k))=0$ except maybe if $k=0,1$. Then $X$
 is regular if and only if $H^1(E(k))=0$ for all $k$, and 
 by Horrock's splitting Criterion, this is equivalent to the 
 splitting of $E$. Thus $X$ is regular if and only if $\varphi_*\Cal
 O_Y$ splits as direct sum of line bundles. \qed
 
 \smallskip
 
 The following proposition tells how the restriction of $\varphi_*\Cal O_X$
 to a smooth curve in 
$|\omega_X|$ splits: 
 
 \proclaim{\hypsplit } Let $W$ be a surface of minimal degree
 $r$, not necessarily 
 smooth, let $X @> \varphi >> W$ be a canonical cover
 of degree $4$
and let
 $C$ be a general smooth irreducible curve  in $|\Cal
 O_W(1)|$.
 If $X$ is regular, 
 then 
 $$(\varphi_*\Cal
 O_X) |_C=\Cal O_{\bold P^1} \oplus  \Cal O_{\bold P^1}(-r-1)  \oplus
 \Cal O_{\bold P^1}(-r-1) \oplus \Cal O_{\bold P^1}(-2r-2) \ .$$   
 \endproclaim
 
 \noindent {\it Proof.} See \GPtrans \ns, Lemma 2.3 for details. 
 \qed
 
 \medskip
 
 Finally we describe more accurately the splitting of $\varphi_*\Cal O_X$ in
 the case $X$ regular:
 
 \proclaim{\splitwist } Let $X @> \varphi >> W$ be a quadruple
 canonical cover. If $\varphi_*\Cal O_X$ splits as a direct sum of
 line bundles, then:  
 
 \item{1)} If $W$ is linear $\bold P^2$, then 
 $$ \varphi_*\Cal O_X= \Cal O_{\bold P^2} \oplus \Cal
 O_{\bold P^2}(-2) \oplus \Cal O_{\bold P^2}(-2) \oplus 
 \Cal
 O_{\bold P^2}(-4) \ ;$$
 
 \item{2)} $W$ is not the Veronese surface; and
 
 \item{3)} if $W$ is a rational normal scroll and $X$ is regular, then
      $$\displaylines{\varphi_*\Cal O_X= \Cal O_W \oplus \Cal
 O_W(-C_0-(m+1)f) \oplus \cr \Cal O_W(-2C_0-(e+1)f) \oplus  
 \Cal O_W(-3C_0-(m+e+2)f) \ ,}$$ where $2m-e$ is the degree of $W$.  
\endproclaim
 
 {\it Proof.} 
 To prove 1) recall that by \regPtwo $X$ is regular. 
 Then by \hypsplit the restriction of $\varphi_*\Cal O_X$ to a line of $W$ is
 
 $$\Cal O_{\bold P^1} \oplus \Cal O_{\bold P^1}(-2) \oplus \Cal O_{\bold
   P^1}(-2) \oplus \Cal O_{\bold P^1}(-4) $$ so 1) is clear.

 For 2) $W$ is isomorphic to $\bold P^2$ and $\Cal O_W(1)= \Cal
 O_{\bold P^2}(2)$.  
 Then by \split $\varphi_*\Cal O_X=\Cal O_{\bold P^2} \oplus E \oplus \Cal
 O_{\bold P^2}(-5)$. Let  
 $C$ be a smooth conic in $W$. If $X$ is irregular, then by \regPtwo
 \ns, $\varphi_*\Cal O_X$ cannot
 split completely. If $X$ is regular, then according
 to \hypsplit
 $$(\varphi_*\Cal O_X)_C=\Cal O_{\bold P^1}  \oplus \Cal O_{\bold P^1}(-5)
 \oplus \Cal O_{\bold P^1}(-5)
 \oplus \Cal O_{\bold P^1}(-10)\ .$$  
 Then if  
 $E=L_1^* \oplus L_2^*$ with  $L_1$ and $L_2$ line bundles, 
 $L_1^*|_C=L_2^*|_C=\Cal O_{\bold P^1}(-5)$, but since $C$ is a conic the
 degrees of $L_1|_C$ and $L_2|_C$ are even integers, so we get a
 contradiction.
 
 For 3) 
 recall that $\omega_W=\Cal O_W(-2C_0-(e+2)f)$ and that
 $\Cal O_W(1)=\Cal O_W(C_0+mf)$. Then it follows by assumption and
 \split  that
 $$\varphi_*\Cal O_X=\Cal O_W \oplus L_1^* \oplus L_2^* \oplus \Cal
 O_W(-3C_0-(m+e+2)f) \ ,$$ and 
 $L_1 \otimes L_2 = \Cal O_W(3C_0+(m+e+2)f)$. Then, if we set
 $L_1=\Cal O_W(a_1C_0+b_1f)$ and $L_2=\Cal O_W(a_2C_0+b_2f)$, it
 follows that
 $$\displaylines{a_1+a_2=3 \cr
 b_1+b_2=m+e+2 \ .}$$
 We show now that $a_i \geq 1$ for $i=1,2$. Since $\varphi$ is induced by
 the complete linear series of $\varphi^*\Cal O_W(C_0+mf)$, then
 $(1-a_i)C_0+(m-b_i)f$ is 
 non effective. Then, if $a_i \leq 0$, $b_i \geq m+1$.     
 On the other hand, since $X$ is regular and $H^1(L_i^*) \subset
 H^1(\varphi_*\Cal O_X)=0$, it follows that $H^1(\Cal
 O_W(-a_iC_0-b_if))=0$. Then, if $a_i \leq 0$, $b_i \leq 1$, hence $m
 \leq 0$, which is impossible, because $m \geq e+1 \geq 1$. Then,
 since $a_1+a_2=3$, $a_i$ is either 
 $1$ or $2$. Let us set $a_1=1$. Then, since $(1-a_1)C_0+(m-b_1)f$
 cannot be effective, $b_1 \geq m+1$. On the other hand, $a_2=2$, and
 since $0=h^1(L_2^*)=h^1(\Cal O_W(-2C_0-b_2f))=h^1(\Cal
 O_W((b_2-e-2)f))=h^1(\Cal O_{\bold P^1}(b_2-e-2))$, then $b_2 \geq
 e+1$. Since $b_1+b_2=m+e+2$, then $b_1=m+1$ and $b_2=e+1$. \qed
 
\medskip 

The purpose of this paper is to study canonical quadruple Galois
covers, and we will focus on them in the next sections. Meanwhile,
\splitwist already yields a fact regarding these covers which is worth
remarking. Note that a canonical Galois cover of the Veronese surface
is flat, 
because the Veronese surface is smooth. Then, being $\varphi$ flat and Galois, 
$\varphi_*\Cal O_X$ splits completely, so we have this

\proclaim{\NoVer } There are no quadruple Galois canonical covers of
the Veronese surface. 
\endproclaim

 \heading 2. Quadruple Galois covers \endheading
 
 In this section we start the study of 
 the case when $\varphi$ is a Galois, quadruple
 cover. We begin by setting up our definitions and recalling, without a
 proof, some useful facts on Galois extensions.

 \proclaim{\defGal } 
 Let $\frak X @> \Phi >> \frak W$ 
 be a finite morphism between normal projective
 varieties $\frak X$ and $\frak W$ and let $G$
 be a finite group 
 acting on $\frak X$ so that $\frak X/G=\frak W$ and 
 $\Phi$ is in fact the projection
 from $\frak X$ to $\frak X/G$. We call $\frak X @> \Phi >> \frak W$ 
 a Galois cover with group
 $G$. 
 \endproclaim

 \proclaim{\Galremtwo } In the previous definition of Galois cover we are not
 assuming $\Phi$ to be flat. However, a quadruple Galois canonical cover $X
 @> \varphi >> W$ of a smooth surface of minimal degree  
is flat, because $W$ is smooth, $X$ is locally Cohen-Macaulay and
 $\varphi$ is finite. 
 \endproclaim

 \proclaim{\Galemma } Let $A \subset B$ be a ring extension of
 integrally closed domains, let $K$ be the field of fractions of $A$ and
 let $L$ be the field of 
 fractions of $B$. Let $G$ be
 a finite group. Then following are equivalent:
 
 \item{1)} $G$ acts on $B$ and $A=B^G$. 
 \item{2)} $K \subset L$ is a normal extension with Galois group $G$
 and $B$ is the integral 
 closure of $A$ in $L$. 
 
 Under the equivalent conditions 1) and 2), $G$ is the group of
 $A$-automorphisms of $B$.
 
 \endproclaim
 
 Now we focus on Galois covers of degree $4$. 
 We recall the algebra structure associated to such
 covers: 
 
 \proclaim{\Galquad } 
 Let $\frak X$ and $\frak Y$ be two algebraic 
 varieties and let $\frak X @> p >> \frak Y$ be a flat,
 Galois cover of degree $4$.

 \item{1)} If $G=\bold Z_4$, then  $p_*\Cal O_\frak X$ splits as
 $$p_*\Cal O_\frak X = \Cal O_\frak Y \oplus L_i^* \oplus L_{-1}^* \oplus
 L_{-i}^*$$ where the line bundles on $\frak Y$,   
 $\Cal O_\frak Y, L_i^*, L_{-1}^*$ and $ L_{-i}^*$ are the eigenspaces of
 $1,i,-1$ and $-i$ respectively.

 There exist effective Cartier divisors $D_{11}, D_{12}, D_{23}$ and $D_{33}$
 on $\frak Y$ such that 
 $D_{11}+D_{23}=D_{12}+D_{33}$ and the following 
 $$\matrix L_i \otimes L_i &=& L_{-1} \otimes \Cal O_Y(D_{11}) \cr
 L_i \otimes L_{-1} & = & L_{-i} \otimes \Cal O_Y(D_{12})\cr
 L_i \otimes L_{-i} & = & \Cal O_Y(D_{11}+D_{23})  
 \cr
 L_{-1}\otimes L_{-1} & = & \Cal O_Y(D_{12}+D_{23})\cr
 L_{-1}\otimes L_{-i} & = & L_i \otimes \Cal O_Y(D_{23})\cr
 L_{-i} \otimes L_{-i} & = & L_{-1} \otimes \Cal O_Y(D_{33})\cr 
 \endmatrix 
 \qquad \Galquadone    
 $$
 and the multiplicative structure of  $p_*\Cal O_\frak X$ is as
 follows: 
 
 $$\displaylines{L_i^* \otimes L_i^* @> \cdot D_{11} >> L_{-1}^* \cr
 L_i^* \otimes L_{-1}^* @> \cdot D_{12} >> L_{-i}^*\cr
 L_i^* \otimes L_{-i}^* @> \cdot D_{11}+D_{23}  >> \Cal O_\frak Y  \cr
 L_{-1}^* \otimes L_{-1}^* @> \cdot D_{12}+D_{23} >> \Cal O_\frak Y \cr
 L_{-1}^*\otimes L_{-i}^* @> \cdot D_{23} >> L_i^*\cr
 L_{-i}^* \otimes L_{-i}^*  @> \cdot D_{33} >> L_{-1}^* \ . \cr}$$

 \smallskip
 
 \item{2)} If $G=\bold Z_2 \times \bold Z_2$, 
 then $p_*\Cal O_\frak X$ splits as
 $$p_*\Cal O_\frak X = \Cal O_\frak Y \oplus L_1^* \oplus L_2^* \oplus
 L_3^* \ ,$$ where $\Cal O_\frak Y$, $L_1^*$,
   $L_2^*$ and $L_3^*$ are eigenspaces
 and there exist effective Cartier divisors 
   $D_1$, $D_2$ and $D_3$ such that 
 $L_i^{\otimes 2}=\Cal O_Y(D_j+D_k)$ and 
 $L_j \otimes L_k = L_i \otimes \Cal O_Y(D_i)$ with $i \neq j$, $j \neq
 k$ and $k 
 \neq i$, and the multiplicative structure of $p_*\Cal O_\frak X$ is as
 follows: 
 $$\displaylines{L_i^* \otimes L_i^* @> \cdot D_j+D_k >> \Cal O_\frak Y \ ,
   \cr
 L_j^* \otimes L_k^* @> \cdot D_i >> L_i^*
 \ .}$$ 
 \endproclaim
 
 \noindent {\it Sketch of proof.} 
The action of $G$ splits $p_*\Cal O_\frak X$ into line bundles.
 Precisely if $G$ is $\bold Z_4$, they  are the eigenspaces of the
 eigenvalues $1,i,-1$ and $-i$ of one generator of the group and if $G$ is
 $\bold Z_2 \times \bold Z_2$ they are the intersection of the eigenspaces
 associated to one of the generators of $G$ with the eigenspaces of
 other generator of $G$. Then it is easy to trace to what eigenspace the
 product of elements of two eigenspaces goes. For instance, if $G$ is
 $\bold Z_4$ it is clear that the product of two elements of $L_i$ is
 an element of $L_{-1}$, the product of an element of $L_i$ and an
 element of $L_{-1}$ is an element of $L_{-i}$, etc. This
   yields the existence of the effective divisors mentioned in the
   statement and the description of the multiplicative structure, the
   relations among the divisors arising from 
  the commutativity and associativity of the product. 
 We refer to \Pa and, for the bi-double case,  
to \text{\Ca \ns;} see also \HM \ns. \qed

 \medskip
 
 We recall also the definition of simple cyclic cover:

 \proclaim{\defnrem} Let $D$ be an effective Cartier
 divisor on a variety 
 $\frak Y$ such that there exists a line bundle $L$ with $L^{\otimes
   n}=\Cal O_\frak Y(D)$. Consider the $\Cal O_\frak Y$-algebra $\Cal O_\frak Y
 \oplus L^{-1} \oplus L^{-2} \oplus \cdots \oplus L^{-n+1}$ with
 multiplication
 $$\displaylines{L^{-m_1} \otimes L^{-m_2} \longrightarrow L^{-(m_1+m_2)} \cr
 \text{ if } \quad 
 0 \leq m_1+m_2 \leq n-1,   \quad \text { and }
 \cr
 L^{-m_1} \otimes L^{-m_2} @> \cdot D >> L^{-(m_1+m_2)+n} \cr \text{
   otherwise. }}$$
 Then we say that $\frak X=$ Spec$(\Cal O_\frak Y
 \oplus L^{-1} 
 \oplus \cdots \oplus L^{-n+1})$ is a
 simple cyclic cover of degree $n$ branched along $D$. 
\smallskip
 A  simple cyclic cover of degree $n$ is Galois and its
 Galois group is $\bold Z_n$. 
 In particular, with the notation of \Galquad \ns, 
 a Galois cover $\frak X @> p >> \frak Y$  with group $\bold Z_4$ is
 simple cyclic branched along $D$ if and only if $D_{11}=D_{12}=0$ and
 $D_{23}=D_{33}=D$ or $D_{11}=D_{12}=D$ and $D_{23}=D_{33}=0$
 \endproclaim

 A quadruple Galois canonical cover over a smooth surface of minimal
 degree 
 has the property that two of the eigenspaces (different from
 $\Cal O_\frak Y$) described in \Galquad have as tensor product the third
 one (see \split \ns.) This is a strong condition that simplifies the 
 algebra structures of the covers. We see what this exactly means in the 
 next two propositions. We start with the case when the  Galois group is $\bold
 Z_2 \times \bold Z_2$. 
 
 \proclaim{\bipro } Let $\frak X$ and $\frak Y$ be algebraic varieties
  and let 
 $\frak X @> p >> \frak Y$ 
 be a flat Galois cover with Galois group $\bold Z_2 \times
 \bold Z_2$. Let $L_1, L_2, L_3, D_1, D_2$ and $D_3$ be as in \text{\Galquad
 \ns,} 2). If $L_1 \otimes L_2 = L_3$, then

\item{1)} $\frak X$ 
 is the fiber product over $\frak Y$ 
 of the flat double covers $\frak X_1 @> p_1 >> \frak Y$ and $\frak X_2 @>
 p_2 >> \frak Y$,
which
 are  branched
 respectively along $D_2$ and $D_1$ and have trace zero module $L_1^*$
 and $L_2^*$ respectively, and $p$ is the natural map from the fiber
 product onto $\frak Y$.
 
 \item{2)} In particular, if $X @> \varphi >> W$ is a canonical 
 cover of a smooth surface of minimal degree 
 and $\varphi$ is Galois with Galois group $\bold Z_2
 \times \bold Z_2$, 
 then $X$ is the fiber product over $W$ 
 of two double covers $X_1 @> p_1 >> W$ and $X_2 @> p_2 >> W$
 and $\varphi$ is the natural map from the
 fiber product to $W$. 
 \endproclaim
 
 {\it Proof.} According to \Galquad
 $$p_*\Cal O_\frak X= \Cal O_\frak Y \oplus L_1^* \oplus L_2^*
 \oplus L_3^*  
 $$ and the $\Cal O_\frak Y$-algebra $p_*\Cal O_\frak X$ has three subalgebras
 $\Cal O_\frak Y \oplus L_1^*$, $\Cal O_\frak Y \oplus L_2^*$ and $\Cal
 O_\frak Y \oplus
 L_3^*$, corresponding to three double covers $\frak X_i @> p_i >>
 \frak Y$ for $i=1,2,3$
 which are branched along $D_j+D_k$ where $i \neq j$, $j \neq k$, $k
 \neq i$. 
 Since $L_1 \otimes L_2=L_3$, then $D_3 =0$, hence $\frak X_1 @> p_1 >>
 \frak Y$ is
 branched along $D_2$ and $\frak X_2 @> p_2 >> \frak Y$ 
 is branched along $D_1$. 
 Then the algebra structure of $p_*\Cal O_\frak X$ 
 described in \Galquad \ns, 2) 
 becomes the tensor product over $\Cal O_\frak Y$
 of the algebras $\Cal O_\frak Y \oplus L_1^*$ and  $\Cal O_\frak Y \oplus
 L_2^*$.

 Now let $X @> \varphi >> W$ be a canonical 
 cover of degree $4$. 
If $\varphi$ is Galois, by \text{\Galquad \ns,}  $\varphi_*\Cal O_X$
 splits as direct sum of line bundles. Then by \split 
 $L_1 \otimes L_2 = L_3$, so the statement follows from the
 first part of 
 \text{\bipro \ns.} \qed

 \proclaim{\Galcyclic} 
  Let $\frak X$ and $\frak Y$ be algebraic varieties and let 
 $\frak X @> p >> \frak Y$ 
 be a flat Galois cover with Galois group $\bold Z_4$ and let $L_i,
 L_{-1}, L_{-i}, D_{11}, D_{12}, D_{23}$ and $D_{33}$ be as in \Galquad \ns. 
 If $L_{\tau(i)}
 \otimes L_{\tau(-1)}= L_{\tau(-i)}$ 
 for some permutation $\tau$ of $\{i,-1,-i\}$,
  then 
 \item{1)} either 
 $L_{i} \otimes L_{-1} = L_{-i}$ (i.e., $D_{12}=0$) or $L_{-1} 
  \otimes L_{-i} = L_i$ (i.e.,  $D_{23}=0$). 

 \item{2)} In the case $L_{i} \otimes L_{-1} = L_{-i}$, the cover $\frak X @> p
 >> \frak Y$ is the composition of a flat double cover $\frak X' @> p_1 >>
 \frak Y$ branched along $D_{23}$ followed by a flat double cover $\frak X
 @> p_2 >> \frak X'$, 
 branched along $p_1^*D_{11}$ and the ramification
 locus of $p_1$. Moreover, the trace zero module of $p_2$ is $p_1^*L_{i}^*$
 and the trace zero module of $p_1$ is $L_{-1}^*$.  
The surface $\frak X'$ is the quotient of $\frak X$ by the
 unique subgroup of index $2$ of $\bold Z_4$. If $L_{-1}
 \otimes L_{-i} = 
 L_i$ there exists an analogous decomposition of $p$.
 
 \item{3)} Assume that $\frak Y$ is locally Gorenstein. If
   $L_{i} \otimes L_{-1} = L_{-i}$, then  $\frak X$ is locally Gorenstein and 
 $\omega_\frak
   X=p^*(\omega_\frak Y
   \otimes L_{-i}).$   If $L_{-1}
 \otimes L_{-i} = 
 L_i$, then  $\frak X$ is locally Gorenstein and $\omega_\frak
   X=p^*(\omega_\frak Y
   \otimes L_{i}).$

 \item{4)} In particular, if $X @> \varphi >> W$ is a canonical 
 cover of degree $4$ of a smooth surface of minimal degree
 and $\varphi$ is Galois with Galois group $\bold Z_4$,  
 then $X @> \varphi >> W$ is a cover satisfying 1), 2) and 3) above and either
 $L_i$ or $L_{-i}$ 
 is isomorphic to $\omega_W^*(1)$.    
 \endproclaim
 
 {\it Proof.} 
 By assumption either
 $$\displaylines{L_i \otimes L_{-1} = L_{-i} \text{ or }\cr
 L_{i} \otimes L_{-i} = L_{-1} \text{ or } \cr
 L_{-1} \otimes L_{-i}= L_{i} \ .}$$
 Since $\frak X$ is normal, $L_{i} \otimes L_{-i} = L_{-1}$ is not
 possible.   
 Then either $L_i \otimes L_{-1} = L_{-i}$ or $L_{-1} \otimes
 L_{-i}= L_{i}$ so in \Galquad
  either $D_{12}=0$ or $D_{23}=0$. 
 If, for example, $D_{12}=0$, then
 looking locally at the
 multiplicative structure of $p_*\Cal O_X$
 yields the description of the trace zero module and branch loci given
 in 2). 
 Now, $L_{-1}^*$ is the trace zero module of $p_1$ and $p_1^*L_i^*$ is the
 trace zero module of $p_2$, so  
 $\omega_\frak X = p^*(\omega_\frak Y \otimes L_{i} \otimes L_{-1}) =
  p^*(\omega_\frak Y \otimes L_{-i})$.
 This proves 3).
Finally,  
 from \split we know also that, if $X @> \varphi >> W$ is a canonical 
 cover of degree $4$, 
 then $L_{\tau(i)} \otimes L_{\tau(-1)} = L_{\tau(-i)}$ for some
 permutation $\tau$ of $\{i,-1,-i\}$, so $X @> \varphi  >> W$ satisfies
 the hypothesis of the statement.  Then by 1) either 
$L_i \otimes L_{-1} = L_{-i}$ or 
$ L_{-1} \otimes L_{-i}= L_{i}$. We argue for instance in case $L_i
 \otimes L_{-1} = L_{-i}$. Since $\omega_X=p^*\Cal O_W(1)$, then 3)
 implies that $L_{-i}$ is numerically equivalent to $\omega_W^*(1)$,
 but since $W$ is either $\bold P^2$ or a Hirzebruch surface
 we have $L_{-i}=\omega_W^*(1)$. If $ L_{-1} \otimes L_{-i}= L_{i}$
 arguing similarly we obtain $L_{i}=\omega_W^*(1)$.
 \qed
 
 \proclaim{\Notarem} If  $X @> \varphi >> W$ is a canonical 
 cover of degree $4$ of a smooth surface of minimal degree 
and $\varphi$ is Galois with group $\bold Z_4$,
 then we rename $\{L_{i}, L_{-1}, L_{-1}\}$ 
 as $\{L_1,L_2,L_3\}$ so that $L_1 \otimes L_2 = L_3$, $L_2=L_{-1}$
 and $L_3=\omega_W^*(1)$. Then there exist effective Cartier divisors
 $D_1$ and $D_2$ so that 
 $$\displaylines{L_1 \otimes L_1 = L_2 \otimes \Cal O_W(D_{1}) \cr
 L_1 \otimes L_2 = L_3 \cr
 L_1 \otimes L_3 = \Cal O_W(D_{1}+D_{2})\cr
 L_2\otimes L_2 = \Cal O_W(D_{2})\cr
 L_2\otimes L_3 = L_1 \otimes \Cal O_W(D_{2})\cr
 L_3 \otimes L_3 = L_2 \otimes \Cal O_W(D_{1}+D_{2}) \ .\cr}$$
 This is achieved if we set $D_1=D_{11}$ and $D_2=D_{23}$ in \Galquadone 
 in case $L_1=L_i$ and one can argue analogously if
 $L_1=L_{-i}$. In particular $\varphi$ is  
 simple cyclic if and only if 
  $D_1=0$.
 \endproclaim
 
 Now we make two useful observations which are a sort of converses to \bipro
 and \Galcyclic \ns: 
 
 \proclaim{\biproconv} 
 Let $\frak X$ and $\frak Y$ be normal algebraic varieties.
 \item{1)} If 
 $\frak X @> p >> \frak Y$ 
 is a morphism such that $p$ is the natural map onto $\frak Y$ 
from the fiber product
 over 
 $\frak Y$ 
 of two flat double covers $\frak X_1 @> p_1 >> \frak Y$ and $\frak X_2 @>
 p_2 >> \frak Y$, then $p$ is a Galois cover with Galois group $\bold
 Z_2 \times \bold Z_2$. 
 \item{2)} If in addition  
 $L_2^*$ and $L_1^*$ are the trace zero 
 modules of $p_2$ and $p_1$ respectively, then 
 $$p_*\Cal O_\frak X=\Cal O_\frak Y \oplus L_1^* \oplus L_2^* \oplus
 (L_1^* \otimes L_2^*) $$ and, if $\frak Y$ is locally Gorenstein,
 then  $\frak X$ is locally Gorenstein
 and $\omega_\frak X= p^*(\omega_\frak Y \otimes L_1 \otimes L_2)$. 
 \endproclaim
 
 {\it Proof.} Choosing suitable open sets $U$ and $V$ of $\frak X$ and
 $\frak Y$ we can describe
 the structure of $\Cal O_\frak Y$-algebra of $\Cal O_\frak X$ in local
 coordinates and 
 see that $\bold Z_2 \oplus \bold Z_2$ acts on $\Cal O_U$ and
 its ring of invariants is 
 $\Cal O_V$. 
 Then, since $\frak X$ and $\frak Y$ are normal, by \Galemma 
 the extension $\Cal K(X)/\Cal K(Y)$ is Galois with Galois
 group $\bold Z_2 \times \bold Z_2$, and since since $\frak X$ and
 $\frak Y$ are normal and $p$ is finite, $\Cal O_\frak X$ is the
 integral closure of $\Cal O_\frak Y$ in $\Cal K(X)$, so again by
 \Galemma the morphism $p$ is a Galois cover with group $\bold Z_2 \times \bold
 Z_2$. The computation of the push-forward of $\Cal O_\frak X$ and the
 computation of 
 the canonical of $\frak X$ is
 straight-forward once we recall that the trace zero module of the 
 double cover $p_1$ is
 $L_1^*$ and using the 
 fact that $p=p_1 \circ p_2'$, where $p_2'$ is another double cover
 with trace zero module $p_1L_2^*$. \qed

 \proclaim{\Galcycliconv} 
 Let $\frak X$ and $\frak Y$ be normal algebraic varieties. 
 \item{1)} If 
 $\frak X @> p >> \frak Y$ 
 is the composition of a flat double cover $\frak X' @> p_1 >>
 \frak Y$ branched along a divisor $D_2$,  
 followed by a flat double cover $\frak X
 @> p_2 >> \frak X'$, 
 branched along  the ramification
 locus of $p_1$ and $p_1^*D_{1}$, where $D_1$ is a divisor on $\frak
 Y$,  
 then $p$ is a Galois cover with Galois group $\bold
 Z_4$. 
 \item{2)} If in addition $L_2^*$ is the trace zero module of $p_1$  and 
 $p_1^*L_1^*$ is the trace zero module of $p_2$, then 
 $$p_*\Cal O_\frak X=\Cal O_\frak Y \oplus L_1^* \oplus L_2^* \oplus
 (L_1^* \otimes L_2^*) $$ and, if 
 $\frak Y$ is locally Gorenstein, then $\frak X$ is
 locally Gorenstein and 
 $\omega_\frak
 X=p^*(\omega_\frak Y
 \otimes L_1 \otimes L_2 ).$
 \endproclaim
 
 {\it Proof.}  
 Let $V$ a smooth open set of $\frak Y$ contained in $\frak Y -
 D_1$. Over $V$ the 
 cover $p$ is a simple cyclic cover, and in particular, a Galois cover
 with Galois group $\bold Z_4$. Arguing as in \biproconv we conclude
 that $\frak X @> p >> \frak Y$ is Galois with Galois group $\bold
 Z_4$. The computation of the push-forward of $\Cal O_\frak X$ and of
 the canonical bundle of $\frak X$ is 
 straight-forward from the fact that $p=p_1 \circ p_2$ and the
 knowledge of the trace zero modules of $p_1$ and $p_2$. 
 \qed

 \medskip 
 
 In the remaining of the paper we will classify quadruple 
Galois canonical covers
 of smooth  surfaces $W$ of minimal degree. A priori one could distinguish
 three cases: $W$ is linear $\bold P^2$, 
$W$ is the Veronese surface and $W$ is a smooth 
 rational normal scroll.
 However, as pointed out in \NoVer \ns, there are no canonical quadruple
 Galois 
 covers of the Veronese surface,  so we
 will have only to study the cases of $W$ being linear $\bold P^2$,
 dealt with in the next theorem, and of $W$ being a
 smooth rational normal scroll, dealt with in Section 3.

 \proclaim{\strPtwo} Let $W$ be linear $\bold P^2$ and let 
 $X @> \varphi >> W$ be a canonical Galois 
 cover of degree $4$. 
 \item{1)} If the Galois group of $\varphi$ is $\bold Z_4$, then  
 $\varphi$ is the composition of two flat double covers $X_1 @> p_1 >> W$
 and $X @> p_2 >> X_1$; the cover $p_1$ is branched along a quartic
                 and the cover $p_2$ is 
                 branched along the ramification of $p_1$ and
                 the pullback by $p_1$ of a conic and its trace
                 zero module is $p_1^*\Cal O_{\bold P^2}(2)$. 
                 
 \item{2)} If the Galois group of $\varphi$ is $\bold Z_2 \times \bold
   Z_2$ then  $X$ is 
 the fiber product over
 $W$ of two double covers of linear $\bold P^2$, 
 each of them branched along a quartic, 
 and $\varphi$  is the  
 natural map from the fiber product to $W$.
 
 \smallskip

 Conversely, let $X @> \varphi >> W$ be 
 a cover of linear $\bold P^2$.

\item{1')} If $\varphi$ is the composition of 
two flat double covers $X_1 @> p_1 >> W$
 and $X @> p_2 >> X_1$ as described in 1) above, then 
 $\varphi$ is a Galois canonical
 cover with group $\bold Z_4$.

\item{2')} If $\varphi$ is the  
 natural map to $W$ from the fiber product over $W$ of two double covers as
 described in 2) above, then $\varphi$ is a Galois canonical
 cover with group 
 $\bold Z_2 \times \bold Z_2$.   

 \endproclaim

 \noindent {\it Proof.} 
 By \Galquad the bundle $\varphi_*\Cal O_X$ splits as a sum of line
 bundles.
 Now 1) and 2) follow from   
 \splitwist \ns, 1), 
 \bipro and \Galcyclic \ns. Now we prove the converse. If $\varphi$ is
 a cover as in 1'), then by \Galcycliconv it is Galois with
 Galois group $\bold Z_4$. Likewise if $\varphi$ is a cover as in 2'),
 then by \biproconv it is Galois with Galois group $\bold Z_2 \times
 \bold Z_2$. On the other hand, 
by \biproconv and \Galcycliconv the canonical of
 $X$ is $\varphi^*\Cal O_{\bold P^2}(1)$, so $X$ is a surface of
 general type and $\omega_X$ is base-point-free. Finally to prove that
 $\varphi$ is indeed the morphism induced by $H^0(\omega_X)$ we compare
 $H^0(\varphi^*\Cal O_{\bold P^2}(1))$ and $H^0(\Cal O_{\bold P^2}(1))$
 and see that they are equal. Indeed, to compute $H^0(\varphi^*\Cal
 O_{\bold P^2}(1))$ we push down $\varphi^*\Cal
 O_{\bold P^2}(1)$ to $W$, and since 
 $$\varphi_*\Cal O_X= \Cal O_{\bold P^2} \oplus \Cal O_{\bold P^2}(-2)
 \oplus \Cal O_{\bold P^2}(-2) \oplus \Cal O_{\bold P^2}(-4) $$
 by \biproconv and \Galcycliconv \ns, we obtain $H^0(\varphi^*\Cal
 O_{\bold P^2}(1))=H^0(\Cal O_{\bold P^2}(1))$ as wished. 
 \qed
 
 \medskip
 
 We describe further the Galois covers appearing in \strPtwo \ns:
 
 \proclaim{\Ptwosing  } Let $W$ be linear $\bold P^2$ and let 
 $X @> \varphi >> W$ be a Galois canonical  
 cover of degree $4$. Then,
 \item{1)}  the push-forward of $\Cal O_X$ is as follows:
 $$\varphi_*\Cal O_X= \Cal O_{\bold P^2} \oplus \Cal O_{\bold P^2}(-2)
 \oplus \Cal O_{\bold P^2}(-2) \oplus \Cal O_{\bold P^2}(-4) \ ;$$
 
 \item{2)} the surface $X$ is regular; and 
 
 \item{3)} if the Galois group of $\varphi$ is $\bold Z_4$, then $X$ is
   singular and the mildest possible set of singularities on $X$
   consists of $8$ points of type $A_1$.
 
 \endproclaim

 {\it Proof.} Claim 1) of the corollary has already been shown in the
 proof of \text{\strPtwo \ns.} Claim 2) is straight-forward form 1), since
 $H^1(\Cal O_X)=H^1(\varphi^*\Cal O_X)$. To prove 3) just observe that
 $X_1 @> p_1 >> W$ is branched along a quartic $D_2$ of $\bold P^2$
 and $p_2$ along the 
 ramification of $p_1$ and $p_1^*D_1$, where $D_1$ is a conic of
 $\bold P^2$. If $D_1$ and $D_2$ are both smooth and meet transversaly then
 $X_1$ is smooth and 
 the branch locus of $p_2$ has $8$ singular points of type $A_1$ so $X$
 is smooth except at $8$ points, which are singularities of type
 $A_1$. \qed

 \medskip
 
 We end the section by remarking the existence of examples of covers
 like those appearing in \strPtwo \ns:

 \proclaim{\Ptwoexamples } Let $W$ be linear $P^2$. 
 \item{1)} There exist 
 canonical covers $X @>
 \varphi >> W$ with Galois group $\bold Z_4$ (that is, covers as in
 \strPtwo \ns, 1)) with $8$ singularities of type $A_1$ as only
 singularities. 
 \item{2)} There exist canonical covers $X @>
 \varphi >> W$ with Galois group $\bold Z_2 \times \bold Z_2$ 
 (that is, covers as in
 \strPtwo \ns, 2))
 with $X$ smooth.
 \endproclaim
 
 {\it Proof.} We first deal with 1). 
 By the converse part in \strPtwo  we just have
 to construct a composition of double covers $X @ > p_2 >>
 X_1$ and $X_1 @> p_1 >> W$ branched along suitable
 divisors. Then, using the same notation of the proof of \Ptwosing and by the
 argument there, it suffices to choose $D_1$ a smooth quartic and $D_2$
 a smooth conic of $\bold P^2$ meeting transversally. This is possible
 by Bertini. Now for 2) again by the converse part in \strPtwo  we just have
 to construct the fiber product $\varphi$ of two double covers $X_1 @> p_1 >>
 \bold P^2$ and $X_2 @> p_2 >> \bold P^2$ branched
 along 
 suitable
 divisors $D_2$ and $D_1$, which are two quartics in this case. 
 Then $\varphi=p_1 \circ p_2'$ where $X @> p_2' >> X_1$ is a double
 cover branched along $p_1^*D_1$.  If we choose $D_1$ and $D_2$
 smooth and meeting tranversally (which again is possible by
 Bertini), then $X_1$ is smooth and so is the branch locus of $p_2'$,
 hence $X$ is also smooth. Note that one can construct examples of $X$
 with worse singularities by allowing $D_1+D_2$ to have worse
 singularities. 
 \qed

 \heading 3.  Galois covers of smooth rational normal scrolls
 \endheading 
  
 At the end of Section 2 we saw that quadruple Galois canonical covers
 of $\bold P^2$ can have Galois group $\bold Z_4$ or $\bold Z_2 \times
 \bold Z_2$, but, in the first case, they could not be simple cyclic. 
 The same thing happens for the covers we study in this section and 
 we start by showing that there do not exist simple cyclic 
 quadruple canonical  
 covers $X$ of smooth rational normal scrolls. 
 This was known if $X$ is regular and $W$ 
 is a surface of minimal degree, whether smooth or singular, as consequence
 of the non-existence of canonical simple cyclic covers of degree bigger than
 $3$ (see \GPtrans \ns). 
 In the next
 theorem  we prove the non existence of simple cyclic covers of
 degree bigger than $3$ when $X$
 is an arbitrary surface of general type and $W$ is a smooth rational
 scroll. 
 
 \proclaim{\scrbid} Let  $W$ be a 
smooth 
rational normal scroll,
 let $X @> \varphi >> W$ be a 
 canonical cover of degree $n$. If $X @> \varphi >> W$ is a
 Galois cover and $n \geq 4$, then $\varphi$ is not simple cyclic. 
 \endproclaim 
 
 {\it Proof.} Let us assume $\varphi$ is simple cyclic. Then 
 $$\varphi_*\Cal O_X=\Cal O_W \oplus L^{-1} \oplus \cdots
 \oplus L^{-n+1} \ .$$
 Recall that $W$ is isomorphic to $\bold F_e$. 
 On the one hand $\omega_X=\varphi^*\Cal O_W(1)=\varphi^*\Cal
 O_W(C_0+mf)$, with 
$m \geq e+1$. 
 On the other hand,
 $\omega_X=\varphi^*(\omega_W \otimes L^{-n+1})=\varphi^*(\Cal
 O_W(-2C_0-(e+2)f) 
 \otimes L^{-n+1})$. Thus $\varphi^*L^{-n+1}=\varphi^*\Cal
 O_W(3C_0+(m+e+2)f)$, so $L^{-n+1}$ and $\Cal
 O_W(3C_0+(m+e+2)f)$ are numerically equivalent in $W$. Since $W$ is a
 Hirzebruch surface, $L^{-n+1}=\Cal
 O_W(3C_0+(m+e+2)f)$ and $3$ and $m+e+2$ are both multiple of $n-1$. 
 Since $n \geq 4$ by assumption, this makes $n=4$. In that case, 
 since $\varphi$ is induced by the complete series of $\Cal
 O_W(C_0+mf)$, then $\Cal O_W((m-\frac 1 3 (m+e+2))f)$ should be non
 effective, i.e., $m-\frac 1 3 (m+e+2) <0$, which is the same as
 $2m-e-2 <0$. Then, since $m \geq e+1$, we get the contradiction $e <
 0$ and $\varphi$  cannot be
 simple cyclic. \qed 
 
 \medskip

 After \scrbid we summarize now the status of the existence of simple
 cyclic canonical covers in the following theorem. 
 To see the scope of the result, we remark that
 \propnoprime implies the non-existence of 
 Galois canonical covers of prime degree $p$ of smooth scrolls, $\bold P^2$ or
 the Veronese surface, if $p \geq 5$. 
If in addition $X$ is regular, next theorem
 assures that, if $p \geq 5$, then  
 there are no  
 Galois canonical covers of prime degree $p$ of any surface of minimal degree.
 
 \proclaim{\propnoprime} Let  $W$ be a surface of minimal degree, not
 necessarily smooth, and let $X @> \varphi >> W$ be 
 a Galois canonical cover. 
 If $X$ is regular or $W$ is smooth, and if $\varphi$ is simple cyclic,
   then 
 deg $\varphi \leq 3$.
\endproclaim

 {\it Proof.}  If $X$ is regular, the result follows from
  \GPtrans \ns, Corollary 3.2. So we will assume that $X$ is
 irregular and $W$ smooth. The surface $W$ cannot be isomorphic to
  $\bold P^2$, for if it were, since $\varphi$ is simple cyclic,
  $\varphi_*\Cal O_X$ would split completely, and so $X$ would be regular.  
 Thus  
 $W$ is a smooth rational normal scroll. Then we conclude that
 deg$\varphi \leq 3$
   by 
 applying \scrbid  \ns. 
 \qed

\proclaim{\cornoprime}  Let  $W$ be a surface of minimal degree, not
 necessarily smooth, and let $X @> \varphi >> W$ be 
 a Galois canonical cover. If $X$ is regular or $W$ is smooth and  
$\varphi$ is a Galois canonical cover of prime degree, 
 then deg $\varphi \leq 3$.  
 \endproclaim
 
{\it Proof.}  If deg$\varphi=p$ is prime, then $G$ is cyclic of order $p$, so
  the stabilizer 
  of any $x \in X$ is either $\{id\}$ or $G$, so $\varphi$ is simple
  cyclic. \qed

 \medskip
 
 These results hint towards a positive solution to the following 
 very interesting question regarding Galois canonical covers of prime 
 degree bigger than $3$.
 
 \proclaim{\ques} If $X @> \varphi >> W$ is a Galois canonical cover of 
 prime degree, is then deg $\varphi \leq 3$? 
 \endproclaim

 \scrbid showed that there does not exist simple cyclic canonical
 covers of degree bigger than $3$ of smooth rational normal scrolls. 
 The construction of non-simple cyclic
 canonical covers of surfaces is not easy and there are not
 many examples, to the best of the authors knowledge. We now prove 
 a classification theorem for non-simple cyclic canonical covers $X @>
 \varphi >> W$ of degree $4$,
 splitting the cases $X$ regular and $X$ irregular.  
 \nosimcyclreg  and  \nosimcyclirreg below and \exampcyc 
 show that they do 
 exist, but, as \nosimcyclsing tells,  
 they are always singular having at best singularities of
 type $A_1$.
 
 \proclaim{\nosimcyclreg }
 Let $X
                 @> \varphi >> W$ be a Galois canonical cover of degree
                 $4$ and Galois group 
                 $\bold Z_4$ with
                  $X$ regular and let $W=S(m-e,m)$ be a smooth rational normal
 scroll. Then  
 $\varphi$ is the composition of two flat double covers
                 $X_1 @> p_1 >> W$  and $X @> p_2 >> X_1$ which are as
                 follows:
 
 \smallskip
  
 The cover $p_1$
                 is branched along a divisor $D_2$ on $W$. 
 
 \smallskip
 
 The cover $p_2$ is 
                 branched along the ramification of $p_1$ and
                 $p_1^*D_1$ and has trace zero module $p_1^*\Cal
                 O_Y(-\frac 1 2 D_1 - \frac 1 4 D_2)$, where $D_1$ is a
                 divisor on $W$. 
 
 \smallskip
 
 The scroll $W$ and the divisors $D_1$ and $D_2$ satisfy:

 \smallskip
 \itemitem{1)} $W=S(m-e,m)$, $0 \leq e \leq 2$, $m \geq
                     e+1$, $D_1 \sim (2m-e+1)f$, 
                    $D_2 \sim 4C_0+(2e+2)f$
  \itemitem{2)} $W=S(1,1)$ (i.e., $W$ is a smooth quadric
                        hypersurface in $\bold 
                        P^3$), $D_1 \sim 3C_0$ and $D_2 \sim 2C_0+4f$
                         
 \medskip

Conversely, let
 $X @> \varphi >> W$ be the
 composition of two flat double covers $X_1 @> p_1 >> W$ and $X @> p_2 >>
X_1$ as described above,  
then $\varphi$ 
                 is a Galois canonical cover of $W$ with Galois group 
                 $\bold Z_4$ and $X$ is regular.    
 
 \endproclaim

{\it Proof.} \Galcyclic \ns, 4) says that 
$\varphi$ is the composition of two double covers
                 $X_1 @> p_1 >> W$ branched along a divisor $D_2$
                 and $X @> p_2 >> X_1$,
                 branched along the ramification of $p_1$ and
                 $p_1^*D_1$ and, according to \Notarem \ns, 
$$\varphi_*\Cal O_X = \Cal O_W \oplus L_1^* \oplus L_2^* \oplus L_3^*
 \ ,$$ 
 where $L_2^*$ is
 the eigenspace for $-1$, $L_1^*$ is the eigenspace for either $i$ or
 $-i$, $L_3^*$ the eigenspace for $-i$ or $i$ and $L_1 \otimes L_2 =
 L_3$. Moreover, $L_3 =\omega_W^*(1)$, $L_1^{\otimes 2}=L_2 \otimes \Cal
 O_W(D_1)$, $L_2^{\otimes 2}=\Cal O_W(D_2)$, $p_1^*L_1^*$ is the
 trace zero module of $p_2$ and $L_2^*$ is the trace zero module of
 $p_1$. From this we obtain that the trace zero module of $p_2$ is
 $p_1^*L_1^*=p_1^*\Cal O_Y(-\frac 1 2 D_1 - \frac 1 4 D_2)$.   
 
 \smallskip
 
 Now we show that $W, D_1$ and $D_2$ satisfy 1) or 2). Recall that $W$
 is isomorphic to $\bold F_e$. Since $X$ is
 regular we can apply \splitwist \ns, 3). Then, since $L_3
 =\omega_W^*(1)$ we have either 
 $L_1^{*}=\Cal O_W(-C_0-(m+1)f), L_2^{*}=\Cal O_W(-2C_0-(e+1)f)$ and $  
 L_3^{*}= \Cal O_W(-3C_0-(m+e+2)f)$ or 
 $L_2^{*}=\Cal
 O_W(-C_0-(m+1)f), L_1^{*}=\Cal O_W(-2C_0-(e+1)f)$ and $  
 L_3^{*}= \Cal O_W(-3C_0-(m+e+2)f)$ 
 
\smallskip
 
 {\it Case 1:} $L_1^{*}=\Cal
 O_W(-C_0-(m+1)f), L_2^{*}=\Cal O_W(-2C_0-(e+1)f)$. 
 From the previous description, $\varphi$
 is the composition of $X_1 @> p_1 >> W$, where $p_1$ is a double cover
 branched along a divisor $D_{2}$ linearly equivalent to $4C_0+(2e+2)f$
  and $X @> p_2 >> X_1$, where $p_2$ is a double
 cover branched along the ramification of $p_1$ and $p_1^*D_{1}$,
 where $D_{1}$ is linearly equivalent to $(2m-e+1)f$.
  Recall that $X$ is normal, hence the components of $D_1+ D_2$ have
 multiplicity $1$  and in particular, the fixed part of
 $|4C_0+(2e+2)f|$ contains $C_0$ with multiplicity at most $1$. Thus $e
 \leq 2$.  
 
 \smallskip
 {\it Case 2:} $L_2^{*}=\Cal
 O_W(-C_0-(m+1)f), L_1^{*}=\Cal O_W(-2C_0-(e+1)f)$.
 Again from the description above, 
 $\varphi$
 is the composition of $X_1 @> p_1 >> W$, where $p_1$ is a double cover
 branched along a divisor $D_{2}$ linearly equivalent to $2C_0+(2m+2)f$
 and $X @> p_2 >> X_1$, where $p_2$ is a double
 cover branched along the ramification of $p_1$ and $p_1^*D_{1}$,
 where $D_{1}$ is linearly equivalent to $3C_0+(2e-m+1)f$.
  Recall that $X$ is normal, hence the components of $D_1+ D_2$ have
 multiplicity $1$  and, in particular, the fixed part of
 $|3C_0+(2e-m+1)f|$ contains $C_0$ with multiplicity at most $1$. Thus 
$m \leq 1$, and in fact $m=1$ and $e=0$. 
 So finally, $W=\bold F_0$, $D_1$ is linearly equivalent to $3C_0$ and
 $D_2$ is linearly equivalent to $2C_0+4f$. 
 
 \smallskip
 
 Now we prove the converse. By \Galcycliconv the morphism
  $\varphi$ is a Galois cover with group $\bold
 Z_4$. \Galcycliconv also tells us 
  that $\omega_X=\varphi^*(\omega_W \otimes L_1 \otimes
 L_2)$. Since  $L_1 \otimes L_2 = \Cal O_W(3C_0+(m+e+2)f) =
 \omega_W^*(1)$ in either 1) or 2), we have $\omega_X=\varphi^*\Cal
  O_W(1)$, so $X$ is a surface of general type with base-point-free
  canonical bundle. 
 Now, again by \Galcycliconv \ns,  in both 1) and 2) we have
 $$\displaylines{\varphi_*\Cal O_X = \Cal O_W \oplus \Cal
   O_W(-C_0-(m+1)f) \oplus \cr
 \Cal O_W(-2C_0-(e+1)f) \oplus \Cal O_W(-3C_0-(m+e+2)f) \ . \quad
\rlap
 \nosimcyclone} $$ 
 Then  to see that $\varphi$ is the canonical
  morphism of $X$ we compare $H^0(\omega_X)$$=H^0(\varphi^*\Cal O_W(1))$ and
  $H^0(\Cal O_W(1))$. The group $H^0(\varphi^*\Cal O_W(1))$ can be
  computed pushing $\varphi^*\Cal O_W(1)$ down to $W$ and using
  \nosimcyclone \ns, and one sees at once that $H^0(\omega_X)=H^0(\Cal
  O_W(1))$. Finally we see that $H^1(\Cal O_X)=0$ also by pushing
  down to $W$ and using \text{\nosimcyclone \ns. \qed}
 
 \proclaim{\nosimcyclirreg } Let $X
                 @> \varphi >> W$ be a canonical Galois cover of degree
                 $4$ and Galois group 
                 $\bold Z_4$ with
                  $X$ irregular and let $W=S(m-e,m)$ be a smooth
                 rational normal 
 scroll. Then the irregularity of $X$ is
 $q(X)=1$ and $W$ is isomorphic to $\bold F_0$. Moreover,  
 $\varphi$ is the composition of two flat double covers
                 $X_1 @> p_1 >> W$  and $X @> p_2 >> X_1$ which are as
                 follows:
 
 \smallskip
  
 The cover $p_1$
                 is branched along a divisor $D_2$ on $W$. 
 
 \smallskip
 
 The cover $p_2$ is 
                 branched along the ramification of $p_1$ and
                 $p_1^*D_1$ and has trace zero module $p_1^*\Cal
                 O_Y(-\frac 1 2 D_1 - \frac 1 4 D_2)$, where $D_1$ is a
                 divisor on $W$. 
 
 \smallskip
 
 The scroll $W$ and the divisors $D_1$ and $D_2$ satisfy:

 \itemitem{1)} $W=S(m,m)$, $m \geq 1$,
 $D_1 \sim (2m+4)f$ and $D_2 \sim 4C_0$.

 \itemitem{2)} $W=S(1,1)$, $D_1 \sim 6C_0$ and $D_2 \sim 4f$.

 \medskip
 
 Conversely, if 
$X @> \varphi >> W$ is the
 composition of two double covers $X_1 @> p_1 >> W$ and 
$X @> p_2 >> X_1$ as described above,  
 then $\varphi$ 
                 is a canonical Galois cover of $W$ with Galois group 
                 $\bold Z_4$ and $X$ is irregular.

 \endproclaim

 {\it Proof.} Using \Galcyclic \ns, 4) and \Notarem as in the proof of
 \nosimcyclreg we conclude that 
 $\varphi$ is the composition of two double covers. The first one is 
                 $X_1 @> p_1 >> W$, is branched along a
                 divisor $D_2$ and has trace zero module $L_2^*$. 
 The second cover  is 
                 $X @> p_2 >> X_1$, is 
                 branched along the ramification of $p_1$ and
                 $p_1^*D_1$ and has trace zero module  $p_1^*L_1$. 
 Then $L_1^{\otimes 2}=L_2 \otimes \Cal
 O_W(D_1)$, $L_2^{\otimes 2}=\Cal O_W(D_2)$ and moreover, $L_1 \otimes
 L_2 = L_3
 =\omega_W^*(1)$.  
 Recall that $W=\bold F_e$ 
and let $L_1=\Cal O_W(a_1C_0+b_1f)$ and $L_2=\Cal O_W(a_2C_0+b_2f)$. 
 Then we have
 $$\displaylines{a_1 + a_2 = 3 \cr
 b_1 + b_2 = m + e + 2 }$$  
 Since $L_1^{\otimes 2} \otimes L_2^* = \Cal
 O_W(D_1)$ and $L_2^{\otimes 2}=\Cal O_W(D_2)$ are  effective then 
 $a_2 \leq 2a_1$, $b_2 \leq 2b_1$ and $a_2, b_2 \geq 0$. Then $a_1, b_1
 \geq 0$ also. Moreover, $a_1, b_1 \geq 1$, otherwise we will
 contradict $a_1 + a_2 = 3$ or $b_1 + b_2 = m + e + 2$. 
 We see that $a_1$ cannot be $2$. If $a_1=2$,  
 then $D_1 \sim 3C_0+(2b_1-b_2)f$, $D_2 \sim 2C_0+2b_2f$, $L_1 = \Cal
 O_W(2C_0+b_1f)$ and $L_2=\Cal O_W(C_0+b_2f)$. Since $X$
 is irregular, and $H^1(\Cal O_W)=H^1(L_2^*)=H^1(L_3^*)=0$, then
 $H^1(L_1^*) \neq 0$. This implies $b_1 \leq e$. Then $b_1 + b_2 = m +
 e + 2$ implies $b_2 \geq m + 2$. On the other hand, since $X$ is
 normal, $C_0$ has at most multiplicity $1$ in the fixed part of
 $|3C_0+(2b_1-b_2)f|$, and this implies $2b_1-b_2-2e \geq 0$. Then we
 have $2e -(m+2) -2e \geq 0$, which is a contradiction. Then the only
 possibilities are $a_1=1$ or $a_1=3$
 \smallskip
 
{\it Case 1:} $a_1=1$. Then $a_2=2$ and $X$ irregular implies $b_2
 \leq e$, since 
 $L_2^*$ has to be special. The fact that $X$ is normal implies $-3e+2b_2
 \geq 0$, since 
 $D_2$ cannot have $2C_0$ as a fixed component. Then $b_2=e=0$, and,
 summarizing, 
 $e=0, a_1=1, a_2=2, b_1=m+2, b_2=0$. This implies $D_1 \sim (2m+4)f$
 and $D_2 \sim 4C_0$. 

 \smallskip
 
{\it Case 2:} $a_1=3$. 
 Then  $D_1 \sim 6C_0 + (2b_1-b_2)f$, and since $X$ is normal, 
$C_0$ has at most multiplicity $1$ in the fixed part of $|D_1|$, hence 
 $2b_1 - b_2 -5e \geq 0$. Now since
 $H^0(\varphi^*\Cal O_W(C_0+mf))=H^0(\Cal O_W(C_0+mf))$, we have that $b_2 >
 m$, hence $b_1 < e+2$. Then we get $2e+4-m-5e=-3e-m+4 > 0$. But $m
 \geq e+1$, so this gives $-4e+3 > 0$, hence $e=0$. In this  case,
 $b_1=1$ and $m \geq 1$. Then $b_2=m+1 \geq 2$, and since $D_1$ is
 effective, $b_2=2$ and $m=1$. Summarizing, 
 $e=0$, $m=1$, $a_1=3$, $a_2=0$, $b_1=1$ and $b_2=2$. 
 Then $D_1 \sim 6C_0$ and $D_2 \sim 4f$. 
 \smallskip

 Finally the irregularity of $X$ is $h^1(\Cal O_X)=h^1(\varphi_*\Cal
 O_X)$. We observe that the
 computation of $a_1, a_2, b_1$ and $b_2$ yields that $\varphi_*\Cal
 O_X$ is 
 $$\displaylines{ \Cal O_W \oplus \Cal O_W(-C_0-(m+2)f)  \oplus
 \Cal O_W(-2C_0) 
 \oplus \Cal O_W(-3C_0-(m+2)f) \text{ or }\cr
 \Cal O_W \oplus \Cal O_W(-3C_0-f)  \oplus
 \Cal O_W(-2f) 
 \oplus \Cal O_W(-3C_0-3f) \quad \rlap \quad \nosimcycltwo }$$
 so $h^1(\Cal O_X)=h^1(\Cal O_W(-2C_0))$ or $h^1(\Cal O_X)=h^1(\Cal
 O_W(-2f))$, and in both cases, equal to $1$. 
 \smallskip
 We prove now the converse. \Galcycliconv implies that $\varphi$ is
 Galois with Galois group $\bold Z_4$. Now if $L_2^*$ is the trace zero
 module of $p_1$ and $p_1^*L_1^*$ is the trace zero module of $p_2$,
 then $L_1 \otimes L_2 = \Cal O_W(\frac 1 2 D_1 + \frac 3 4 D_2)$. Then
 if $D_1$ and $D_2$ are as in 1) or 2), $L_1 \otimes L_2
 =\omega_W^*(1)$. Then  \Galcycliconv implies also that
 $\omega_X=\varphi^*\Cal O_W(1)$, therefore $X$ is a surface of general
 type with base-point-free canonical bundle. Finally \Galcycliconv
 tells us also that $\varphi_*\Cal O_X$ is as in \nosimcycltwo so
 arguing like in the end of the proof of \nosimcyclreg we see that
 $\varphi$ is the canonical morphism and $X$ is irregular. 
 \qed

 \proclaim{\nosimcyclsing } Let $W$ be a  smooth rational scroll of
 degree $r$ and let $X @> \varphi >> W$ be a Galois canonical
 cover with Galois group $\bold Z_4$ (i.e.,
 a cover like the ones classified in \nosimcyclreg and \nosimcyclirreg
 \ns). Then $X$ is singular. Moreover,

\item{1)} if $X$ is regular, then 
the mildest possible set of
 singularities on $X$ 
 consists of $4(r+1)$ singular points of type $A_1$ and

\item{2)}  if $X$ is irregular
 the singularities of $X$ are exactly $4(r+4)$ points of type $A_1$.   
 \endproclaim
 
 {\it Proof.} The proof goes as the proof of \Ptwosing \ns, 3). In this
 ocassion, $D_1 \cdot D_2 = 4(r+1)$ if $X$ is regular (see
 \nosimcyclreg \ns, 1) and 2)) and  $D_1 \cdot D_2 = 4(r+4)$ (see
 \nosimcyclirreg \ns, 1) and 2)). $X$ has the mildest possible set of
 singularities if $X_1$ is smooth and the branch locus of $X @> p_2 >>
 X_1$ has the 
 mildest possible set of 
 singularities. This happens if  $D_1$ and $D_2$ are smooth and
 meet tranversally. In this case the branch locus of $p_2$ has only
 singularities of type $A_1$, and so does $X$. Now, if $X$ is
 irregular \nosimcyclirreg together with the fact that $X$ is normal
 implies that $D_1$ is a union of distinct lines of one of the
 fibrations of $\bold P^1 \times \bold P^1$ and $D_2$ is a union of
 distinct lines of the other 
 fibration, so $D_1$ and $D_2$ are smooth and meet transversally in
 any case. \qed

 \medskip

 Now we proceed to classify quadruple Galois canonical 
 covers of smooth rational
 normal scrolls with Galois group $\bold Z_2 \times \bold Z_2$. 
 Having in account \bipro
 we already know
 that 
  they are the fiber product of two double covers. Thus to
 complete their description we will find out what the branch
 loci of the double covers are.    
 We start with the case where $X$ is regular:  
 
 \proclaim{\strsmscreg } 
 Let $W=S(m-e,m)$ be a
 smooth rational normal scroll.
 If  $X$ is  regular and 
 $X @> \varphi >> W$ is a Galois canonical cover 
 with Galois group $\bold Z_2 \times \bold Z_2$, 
 then $X$ is 
 the fiber product  over
 $W$ of two double covers of $X_1 @> p_1 >> W$ and $X_2 @> p_2 >> W$
 and $\varphi$  is the natural map from  the fiber product to $W$. 
 Let the branch divisors 
 $D_2$, $D_1$  of $p_1, p_2$ 
 be linearly
 equivalent to
 $2a_2C_0+2b_2f$ and $2a_1C_0+2b_1f$ respectively.
Then
$0 \leq e \leq 2$, $m \geq e+1$,
  $a_1=1,  a_2=2, b_1=m+1$ and $b_2=e+1$.
 
 \medskip
 
 Conversely, let $W=S(m,m-e)$ be such that $0 \leq e \leq 2$ and $m
 \geq e+1$ and let $X @> \varphi >> W$ be the natural map to $W$ from 
 the fiber product over $W$ of two flat double covers $p_1$ and $p_2$
 with branch divisors as
 described above.  
 Then $X$ is regular
 and $X @> \varphi 
 >> W$ is a Galois canonical cover with Galois group
 $\bold Z_2 \times \bold Z_2$. 
 \endproclaim
 
 {\it Proof.}
 Since $X$ is regular, \splitwist \ns, 3) yields
 $$\displaylines{\varphi_*\Cal O_X= \Cal O_W \oplus \Cal
 O_W(-C_0-(m+1)f) \oplus \cr \Cal O_W(-2C_0-(e+1)f) \oplus  
 \Cal O_W(-3C_0-(m+e+2)f) \ .}$$
 Then \bipro tells us that $X$ is the fiber product  over
 $W$ of two double covers of $X_1 @> p_1 >> W$ and $X_2 @> p_2 >> W$
 with trace zero modules $L_1^*=\Cal O_W(-2C_0-(e+1)f)$ and $L_2^*=\Cal
 O_W(-C_0-(m+1)f)$ respectively, or equivalently, branched along
 divisors $D_2 \sim 4C_0+2(e+1)f$ and $D_1 \sim 2C_0+2(m+1)f$
 respectively. Thus $a_1=1, a_2=2, b_1=m+1$ and $b_2=e+1$. 
 Recall that $W$ is isomorphic to the Hirzebruch surface
 $\bold F_e$. Since $W$ is smooth, $m \geq e+1$, hence the only thing
 left to prove is $e \leq 2$.  
The covers $p_1$ and $p_2$ fit in the
 commutative diagram
 $$ \matrix X & @> p_1' >> &  X_2 \cr
         @VV p_2' V &  \hskip -.8 truecm @VV p_2 V \cr
  X_1 & @> p_1 >>& W \cr 
 \endmatrix \ $$ 
 where $p_1'$ and $p_2'$ are also double covers. 
 Moreover the branch divisor of $p_1'$ is $p_2^*D_2$. 
 Suppose that $e \geq 3$. Then, since  $D_2$ is
 linearly equivalent to
 $4C_0+(2e+2)f$, $D_2$ has $2C_0$ as a fixed component. Then the
branch divisor of $p_1'$ is nonreduced, so $X$ is non-normal
 and we get a contradiction.  Therefore $e=0,1$ or $2$. 
  
 \medskip
 
 To prove the converse assume now that $X @> \varphi >> W$ is the
 natural map from
 the fiber product over a smooth scroll $W=S(m,m-e)$ of two double covers
 $p_1$ and $p_2$
 of $W$, branched respectively along divisors $D_2$ linearly equivalent
 to $2a_2C_0+2b_2f$ and $D_1$ linearly equivalent to
 $2a_1C_0+2b_1f$. Assume in addition that $0 \leq e \leq 2$, $a_1=1,
 a_2=2, b_1=m+1$ and $b_2=e+1$. 
 Then 
by \biproconv $\varphi$ is a 
 Galois cover with Galois group $\bold Z_2 \times \bold Z_2$ and
 $$\displaylines{\varphi_*\Cal O_X= 
 \Cal O_W \oplus \Cal O_W(-a_1C_0-b_1f)
 \oplus \cr
 \Cal O_W(-a_2C_0-b_2f) \oplus \Cal O_W(-(a_1+a_2)C_0-(b_1+b_2)f) \ .} $$
 A standard
 computation shows that none of the four direct summands of $p_*\Cal
 O_X$ have intermediate cohomology, hence $H^1(\Cal O_X)=0$.
 On the other hand if 
$L_2=\Cal O_W(a_1C_0+b_1f)$ and $L_1=\Cal O_W(a_2C_0+b_2f)$,
 $L_1 \otimes L_2=\Cal O_W(3C_0+(m+e+2)f)=\omega_W^*(1)$. Then by
 \biproconv  
 $\omega_X=\varphi^*\Cal O_W(1)$ so $X$ is a surface of general type
 whose canonical bundle is base-point-free.   
 The only thing left to be shown is that $X @> \varphi >> W$ is the
 canonical morphism of $X$. For that it is enough to see that 
 $H^0(\Cal
 O_W(1))=H^0(\omega_X)$. But
 $$\displaylines{H^0(\omega_X)=H^0(\varphi^*\Cal O_W(C_0+mf))=H^0(
 \Cal O_W(C_0+mf)) \oplus \cr H^0(\Cal O_W((1-a_1)C_0+(m-b_1)f)) 
 \oplus \
  H^0(\Cal
 O_W((1-a_2)C_0+(m-b_2)f)) \oplus \cr 
 H^0(\Cal
 O_W(1-a_1-a_2)C_0+(m-b_1-b_2)f)) \ .}$$ 
 Now, because of the restrictions on $a_1, a_2, b_1$ and $b_2$, $H^0(\Cal
 O_W((1-a_1)C_0+(m-b_1)f))$, 
 $H^0(\Cal
 O_W((1-a_2)C_0+(m-b_2)f))$ and $H^0(\Cal
 O_W(1-a_1-a_2)C_0+(m-b_1-b_2)f)$ vanish. \qed

 \medskip
 
 Now we go on to classify Galois quadruple covers with group $\bold Z_2
 \times \bold Z_2$ when $X$ is
 irregular:

 \proclaim{\strsmscrirr} 
  Let $W$ be a
 smooth rational normal scroll $S(m-e,m)$. 
If $X$ is  
 irregular and
 $X @> \varphi >> W$ is a Galois canonical cover 
 with Galois group $\bold Z_2 \times \bold Z_2$, then $X$ is 
 the fiber product over
 $W$ of two double covers of $X_1 @> p_1 >> W$ and $X_2 @> p_2 >> W$
 and $\varphi$  is the  
 natural map from  the fiber product to $W$. 
 Let the branch divisors 
 $D_2$, $D_1$  of $p_1, p_2$ 
 be linearly
 equivalent to
 $2a_2C_0+2b_2f$ and $2a_1C_0+2b_1f$ respectively.
 Then $e=0$, $m \geq 1$ and
one of the following happens:
 
 \smallskip
 
 \item{1)} $a_1=0$, $a_2=3$, $b_1=m+1$, $b_2=1$. 
 
 \item{2)} $a_1=0$, $a_2=3$, $b_1=m+2$, $b_2=0$. 
 
 \item{3)} $a_1=1$, $a_2=2$, $b_1=m+2$, $b_2=0$. 

\smallskip

In addition, in case 1),  $q(X)=m$; in case 2), $q(X)=m+3$; and in
case 3),  $q(X)=1$.

 \medskip
 
 Conversely, let $X @> \varphi >> W$ be the natural map to $W=S(m,m)$ from 
the fiber product over $W$ of two flat double covers $p_1$ and $p_2$
with branch divisors satisfying 1), 2) or  3) above. 
Then $X$ is irregular and $X @> \varphi
 >> W$ is a Galois canonical cover with Galois group
 $\bold Z_2 \times \bold Z_2$.

 \endproclaim

 {\it Proof.} From \split and \bipro it follows that $X$ is the fiber
 product of two double covers branched along divisors $D_2 \sim
 2(a_2C_0+b_2f)$ and $D_1 \sim 2(a_1C_0+b_1f)$ respectively, that 
 $$\displaylines{\varphi_*\Cal O_X= \Cal O_W \oplus \Cal O_W(-a_1C_0-b_1f)
   \oplus 
 \cr \Cal
 O_W(-a_2C_0-b_2f) \oplus \Cal O_W(-(a_1+a_2)C_0-(b_1+b_2)f) }$$
 and that 
 $\omega_W(-1)=\Cal O_W(-(a_1+a_2)C_0-(b_1+b_2)f)$. Since
 $\omega_W=\Cal O_W(-2C_0-(e+2)f)$, we obtain
 $$\displaylines{a_1+a_2=3 \cr
 b_1+b_2=m+e+2 \ \quad \rlap \strsmscrirrnumber }$$ 
 On the other hand since $D_i$ is effective and linearly
 equivalent to $2(a_iC_0+b_if)$, then $a_i, b_i \geq 0$.
 We set $a_1=0,1$ (in which case, $a_2=3,2$). Since $\varphi$ is
 induced by the complete linear series $\varphi^*\Cal O_W(C_0+mf)$, then
 $H^0(\Cal O_W((1-a_1)C_0+(m-b_1)f))=0$, hence $b_1 \geq m+1$ and from
 \strsmscrirrnumber \ns, 
 $b_2 \leq e+1$. Since both $b_1$ and $b_2$ are nonnegative, $m+1
 \leq b_1 \leq m+e+2$ and $0 \leq b_2 \leq e+1$. 
 
 Now assume $a_1=0$. Then $D_2$ is linearly equivalent to
 $6C_0+2b_2f$. Assume also that $e \geq 1$. Then $2C_0$ is a fixed component of
 $D_2$, therefore by the argument made in the proof of \strsmscreg \ns, $X$
 would be nonnormal, hence, if $a_1=0$, then $e=0$. Then, we have two
 possibilities: first, $b_1=m+1$ and $b_2=1$, and second, $b_1=m+2$ and
 $b_2=0$. In the first case, $q(X)=m$. In the second case,
 $q(X)=m+3$. 
 
 Now assume $a_1=1$. Then $D_2$ is linearly equivalent to
 $4C_0+2b_2f$, and as we argued in the proof of \strsmscreg \ns, if $e \geq
 3$, $X$ would be nonnormal, hence $0 \leq e \leq 2$. Moreover, if $b_2
 < \frac{3e}{2}$, $D_2$ would have $2C_0$ as fixed component and $X$
 would be nonnormal, hence $b_2 \geq \frac{3e}{2}$. 
 Recall also that $b_2 \leq e+1$. Let us now assume that $b_2=e+1$. 
 Then $b_1=m+1$ and in  
 that case
 $$\displaylines{\varphi_*\Cal O_X = \Cal O_W \oplus \Cal O_W(-C_0-(m+1)f)
    \oplus \cr \Cal
 O_W(-2C_0-(e+1)f) \oplus \Cal O_W(-3C_0-(m+e+2)f) \ .}$$
 But then $H^1(\Cal O_W)$, $H^1(\Cal O_W(-C_0-(m+1)f))$, $H^1(\Cal
 O_W(-2C_0-(e+1)f))$ and $H^1(\Cal O_W(-3C_0-(m+e+2)f))$ all vanish,
 hence $X$ would be regular. Therefore $\frac{3e}{2} \leq b_2 \leq
 e$. This implies that $e=0$ and $b_2=0$, in which case
 $b_1=m+2$. Then $q(X)=1$. With this we prove that $W=S(m,m)$ and that
 the only possibilities for $a_i$s, $b_i$s 
are
 1), 2) and 3).
 
 \medskip
 
 To prove the converse, assume now that $X @> \varphi >> W$ 
 is the fiber product over $W$ of two double covers
 $p_1$ and $p_2$
 of $W$, branched respectively along divisors $D_2$ linearly equivalent
 to $2a_2C_0+2b_2f$ and $D_1$ linearly equivalent to $2a_1C_0+2b_1f$
 satisfying one of the cases 1), 2) or 3). 
 Then, by \biproconv \ns, $\varphi$ is a
 Galois cover with Galois group $\bold Z_2 \times \bold Z_2$ and
 $$\displaylines{\varphi_*\Cal O_X= 
 \Cal O_W \oplus \Cal O_W(-a_1C_0-b_1f)
 \oplus \cr
 \Cal
 O_W(-a_2C_0-b_2f) \oplus \Cal O_W(-(a_1+a_2)C_0-(b_1+b_2)f)} $$
 with $a_1,a_2, b_1, b_2$ satisfying 1), 2) or 3). Computing the cohomology 
 of the four direct summands of $\varphi_*\Cal
 O_X$ in each case shows that $X$ is always irregular.
 We are going to see that
 $\omega_X=\varphi^*\Cal O_W(1)$. This follows from \biproconv as in
 the proof of 
 \strsmscreg  \ns. Then in
 particular $\omega_X$ is base-point-free and $X$ is a surface of
 general type. 
 The only thing left to be shown is that $X @> \varphi >> W$ is the
 canonical morphism of $X$. For that it is enough to see that 
 $H^0(\Cal
 O_W(1))=H^0(\omega_X)$. This follows because
 $H^0(\omega_X)=H^0(\varphi_*\varphi^*\Cal O_W(1))$, which is equal to
 $H^0(\Cal O_W(1))$ by the same computation as in the proof of
 \text{\strsmscreg \ns.} \qed 
 
 \medskip
 In the next corollary of \classcroll (see also \splitwist \ns, 3) if $X$ is
 regular) we summarize the splitting of
 $\varphi_*\Cal O_X$ if $X @> \varphi >> W$ is a Galois canonical cover
 of degree $4$ of a smooth rational normal scroll: 

 \proclaim{\splitscr } Let $X @> \varphi >> W$ be a canonical Galois cover
 of degree $4$, with Galois group $G$, of a smooth rational
 normal scroll isomorphic to $\bold 
 F_e$ and embedded by
 $|C_0+mf|$. 
 
 \item{1)} If $X$ is regular, then $e=0,1,2$ and 
  $$\displaylines{\varphi_*\Cal O_X= \Cal O_W \oplus \Cal
 O_W(-C_0-(m+1)f) \oplus \cr \Cal O_W(-2C_0-(e+1)f) \oplus  
 \Cal O_W(-3C_0-(m+e+2)f) \ .}$$
 
 \item{2)} If $X$ is irregular, then $e=0$ and 
 
 \smallskip
 
 \itemitem{2.1)} 
 $\varphi_*\Cal O_X= \Cal O_W \oplus \Cal
 O_W(-C_0-(m+2)f) \oplus  \Cal O_W(-2C_0) \oplus  $ \newline 
 $ \Cal O_W(-3C_0-(m+2)f) \ ;\text{ or }$
 
 \smallskip
 
 \itemitem{2.2)}  
  $\varphi_*\Cal O_X= \Cal O_W \oplus \Cal
 O_W(-(m+1)f) \oplus  \Cal O_W(-3C_0-f) \oplus $ \newline $  
 \Cal O_W(-3C_0-(m+2)f) \ ,\text{with $m=1$ if $G=\bold Z_4$; or}$
 
 \smallskip
 
 \itemitem{2.3)}
 $\varphi_*\Cal O_X= \Cal O_W \oplus \Cal
 O_W(-(m+2)f) \oplus  \Cal O_W(-3C_0) \oplus  
 \Cal O_W(-3C_0-(m+2)f),$ \newline and $G=\bold Z_2 \times \bold
   Z_2$.
 \endproclaim
 
 \medskip

 We end the section showing the existence of canonical covers like the
 ones classified in \classcroll \ns:  
 
 \proclaim{\exampcyc} There exist families of quadruple Galois
 canonical 
covers as in \nosimcyclreg  and \nosimcyclirreg which have
 singularities as mild as possible \text{(see \nosimcyclsing \ns).}
 \endproclaim
 
 {\it Proof.} According to the converse part in \nosimcyclreg and
 \nosimcyclirreg 
 we just have
 to construct a composition of double covers $X @ > p_2 >>
 X_1$ and $X_1 @> p_1 >> W$ branched along suitable
 divisors. Let $D_1$ and $D_2$ be as in \nosimcyclreg  and
 \nosimcyclirreg \ns. 
 By the arguments of \Ptwosing and \Ptwoexamples \ns,  in order for
 $X$ to have singularities as mild as possible and, in any case, only
 $A_1$ singularities, it suffices
 to choose $D_1$ and $D_2$ smooth and meeting transversally. We see
 that such choice is indeed possible. 
 For a cover as in \nosimcyclreg \ns,  1),
 $D_1 \sim (2m-e+1)f$ and $D_2 \sim 4C_0+(2e+2)f$.  Then we choose 
 $D_1$ as the union of $2m-e+1$ different fibers. The
 divisor $4C_0+(2e+2)f$ is base-point-free if $e=0,1$ and if $e=2$ is
 $(3C_0+6f)+C_0$, with $3C_0+6f$ base-point-free and $(3C_0+6f)\cdot
 C_0=0$. Thus by Bertini $D_1$ and $D_2$ can be chosen smooth and
 intersecting transversally. 
 For a cover like \nosimcyclreg \ns, 2),
 both $D_1 \sim 3C_0$ and $D_2 \sim 2C_0+4f$ are base-point-free
 (recall that in this case $W$ is isomorphic to $\bold F_0$),  so 
 we can apply 
 Bertini as before. Finally, in \nosimcyclirreg  
 $D_1$ and $D_2$ are the union of distinct smooth lines
 belonging to the two fibrations of the ruled surface $\bold F_0$. 
 Note that, if $X$ is regular, 
one can construct $X$ with worse singularities allowing $D_1
 + D_2$ to have worse singularities. 
 \qed
 
 \proclaim{\exampbid}
 There exist families of canonical Galois quadruple covers $X @>
 \varphi >> W$ as in  
 \strsmscreg and  \strsmscrirr \ns, 1),  2) and
 3) with $X$ smooth. 
 \endproclaim
 
 {\it Proof.}
 Families satisfying \strsmscreg for $W$ isomorphic
 to $\bold F_0$ and $\bold F_1$ 
 have been constructed in \GPtrans \ns,  Examples 3.4 and 3.5 (see also
 \Pe for an example of a bidouble cover of $\bold P^1 \times \bold
 P^1$) \ns. 
 To construct the remaining examples we argue as in \Ptwoexamples
 \ns. 
By the converse part of \strsmscreg
 and \strsmscrirr we just need to construct the fiber
 product of two double covers $X_1 @> p_1 >> W$ and $X_2 @> p_2 >> W$
 branched along suitable divisors $D_2$ and $D_1$ 
 satisfying  the conditions in the statement of \strsmscreg and 
 \strsmscrirr \ns. 
Precisely if we choose $D_1$ and $D_2$ smooth and meeting
 transversally, $X$ will be smooth. This can be achieved using Bertini
 once we study how 
 the divisors $D_1$ and $D_2$ are in each case.
 Indeed, if we are in the situation of \strsmscreg when $W$ is
 isomorphic to $\bold F_2$, then 
 $D_1 \sim 2C_0+2(m+1)f$ is very ample  and $D_2$, since it is linearly
 equivalent to $4C_0+6f$ is of the form $C_0+D_2'$, with $D_2' \cdot
 C_0=0$ and $D_2'$ base-point-free. 
 
 \smallskip
 
 In case \strsmscrirr \ns, 1), 2) and 3) recall that $W$ is isomorphic
 to $\bold F_0$. In case \strsmscrirr \ns, 1)
 $D_1 \sim 2(m+1)f$, hence it can be chosen as the union  of $2(m+1)f$
 distinct lines in one 
 of the two fibrations of $\bold F_0$
 and $D_2$ is  linearly equivalent to 
 $6C_0+2f$, which  is very ample. 
  \smallskip

 Finally in case \strsmscrirr \ns, 3), 
  $D_1 \sim 2C_0+2(m+2)f$  is very ample
  and $D_2 \sim 4C_0$  can be chosen as the union of $4$ distinct lines of one
  of the fibrations of $\bold F_0$. 
 Note that one can construct $X$ with worse singularities if $D_1 +
  D_2$ is allowed to have worse singularities. \qed

 \heading References \endheading
 
 \item{\Be} A. Beauville, {\it L'application canonique pour les surfaces
 de type general}, Inventiones Math. {\bf 55} (1979), 121-140.
 
 \item{\BS} M.C. Beltrametti and T. Szemberg, {\it On higher order
 embeddings of Calabi-Yau threefolds}, Arch. Math. (Basel) {\bf 74}
 (2000), 221-225.
 
 \item{\Ca} F. Catanese, On the moduli spaces of surfaces of general type. 
 J. Differential Geom. {\bf 19} (1984), no. 2, 483--515.

 \item{\GPcy} F.J. Gallego and B.P. Purnaprajna, {\it Very ampleness and higher
 syzygies for Calabi-Yau threefolds}, Math. Ann. {\bf 312} (1998),
 133--149.
 
 \item{\GPtrans} F.J. Gallego and B.P. Purnaprajna, {\it On the canonical ring
 of covers of surfaces of minimal degree}, 
 Trans. Amer. Math. Soc. {\bf 355} (2003), 2715-2732.

\item{\CR } F.J. Gallego and B.P. Purnaprajna, {\it Classification of
    quadruple canonical covers: Galois case}, 
to appear in C. R. Math. Acad. Sci. Soc. R. Can.

\item{\GPsing} 
F.J. Gallego and B.P. Purnaprajna, {\it Classification
    of quadruple Galois canonical covers, II}, preprint. 
 
 \item{\GPring} F.J. Gallego and B.P. Purnaprajna, {\it Ring generation
 for cyclic and quadruple canonical covers}, in preparation.
 
 \item{\Ga} X. Gang, {\it Algebraic surfaces with high canonical degree},
 Math. Ann, {\bf 274} (1986), 473-483.
 
 \item{\Gr} M.L. Green, {\it The canonical ring of a variety of
 general type}, Duke Math. J. {\bf 49} (1982), 1087--1113.
 
 \item{\HM} D. Hahn and R. Miranda, {\it Quadruple covers of algebraic
 varieties}, J. Algebraic Geom. {\bf 8} (1999), 1--30.
 
 \item{\Hoone} E. Horikawa, {\it Algebraic surfaces of general type with
 small $c^2_1$, I}, Ann. of Math. (2) {\bf 104} (1976), 357--387.
 
 \item{\Hotwo} E. Horikawa, {\it Algebraic surfaces of general type with
 small $c\sp{2}\sb{1}$, II}, Invent. Math. {\bf 37} (1976), 121--155.
 
 \item{\Hothree} E. Horikawa, {\it Algebraic surfaces of general type with
 small $c\sp{2}\sb{1}$, III}, Invent. Math. {\bf 47} (1978), 209--248.
 
 \item{\Hofour} E. Horikawa, {\it Algebraic surfaces of general
 type with small $c\sp{2}\sb{1}$, IV}, Invent. Math. {\bf 50} 
 (1978/79),  103--128.
 
 \item{\Kon} K. Konno, {\it Algebraic surfaces of general type with $c_1^2
 =3p_g-6$}, Math. Ann. {\bf 290} (1991), 77--107.
 
 \item{\MP} M. Mendes Lopes and R. Pardini, {\it Triple canonical surfaces
 of minimal degree}, International J.  Math. {\bf 11} (2000), 553--578.

 \item{\OP} K. Oguiso and T. Peternell, {\it On polarized canonical
 Calabi-Yau threefolds}, Math. Ann. {\bf 301} (1995), 237--248.
 
\item{\Pa} R. Pardini, {\it Abelian covers of algebraic varieties},
  J. Reine Angew. Math.  417  (1991), 191--213.  

 \item{\Pe} U. Persson, {\it Double coverings and surfaces of general
 type}, Algebraic Geometry. (Lect. Notes Math., vol 687. 168-195) {\bf 687}
 (1978).

 \item{\Pu} B.P. Purnaprajna, {\it Geometry of canonical covers with 
 applications to Calabi-Yau threefolds}, in Vector Bundles and 
 Representation Theory, Eds. D. Cutkosky et al., 
Contemporary Mathematics Series AMS {\bf
   322} 
(2003), 107--124.
 
 \item{\T} S. Tan, {\it Surfaces whose canonical maps are of odd degrees},
 Math. Annalen, {\bf 292} (1992), 13-29.

 \enddocument